\documentclass[10pt]{article}
\usepackage[nohead,margin=1.0in]{geometry}
\usepackage{amssymb, amsmath, amsthm,amsfonts}
\usepackage{graphicx,epsfig}
\usepackage[tight]{subfigure}
\usepackage{cite}
\usepackage{times}
\usepackage{float}
\graphicspath{{./pics/}}
\usepackage{color}
\usepackage{epstopdf}
\usepackage{mathrsfs}
\usepackage{threeparttable}
\usepackage{algorithm,algorithmic}
\usepackage{booktabs}
\usepackage{verbatim}
\usepackage[colorlinks,linktocpage,linkcolor=blue]{hyperref}

\numberwithin{equation}{section}
\newtheorem{theorem}{Theorem}[section]
\newtheorem{remark}{Remark}[section]
\newtheorem{proposition}{Proposition}[section]
\newtheorem{lemma}{Lemma}[section]

\newtheorem{assumption}[theorem]{Assumption}
\newtheorem{example}{Example}[section]
\allowdisplaybreaks

\def\II{(\Omega)}

\title{Numerical Recovery of the Diffusion Coefficient in Diffusion Equations from Terminal Measurement\thanks{The work of B. Jin is supported by UK EPSRC grant EP/V026259/1, Hong Kong Research Council grant (No. 14306423) and a start-up fund from The Chinese University of Hong Kong, X. Lu by the National Science Foundation of China (No. 11871385), and that of Z. Zhou by Hong Kong Research Grants Council grant (No. 15303021)
and an internal grant of The Hong Kong Polytechnic University (Project ID: P0045708, Work Programme: 4-ZZP7).}}

\author{Bangti Jin\thanks{Department of Mathematics, The Chinese University of Hong Kong, Shatin, N.T., Hong Kong, P. R. China.
(\texttt{b.jin@cuhk.edu.hk, bangti.jin@gmail.com})} \and Xiliang Lu\thanks{School of Mathematics and Statistics, and Hubei Key Laboratory of Computational Science, Wuhan University, Wuhan 430072, P. R. China (\texttt{xllv.math@whu.edu.cn})} \and Qimeng Quan\thanks{School of Mathematics and Statistics, Wuhan University, Wuhan 430072, P. R. China (\texttt{quanqm@whu.edu.cn})}  \and Zhi Zhou\thanks{Department of Applied Mathematics,
The Hong Kong Polytechnic University, Kowloon, Hong Kong, P. R. China (\texttt{zhizhou@polyu.edu.hk})}}

\begin{document}

\maketitle

\begin{abstract}
In this work, we investigate a numerical procedure for recovering a space-dependent diffusion coefficient in a (sub)diffusion model from the given terminal data, and provide a rigorous numerical analysis of the procedure. By exploiting decay behavior of the observation in time, we establish a novel H{\"o}lder type stability estimate for a large terminal time $T$. This is achieved by novel decay estimates of the (fractional) time derivative of the solution. To numerically recover the diffusion coefficient, we employ the standard output least-squares formulation with an $H^1(\Omega)$-seminorm penalty, and discretize the regularized problem by the Galerkin finite element method with continuous piecewise linear finite elements in space and backward Euler convolution quadrature in time. Further, we provide an error analysis of discrete approximations, and prove a convergence rate that matches the stability estimate. The derived $L^2(\Omega)$ error bound depends explicitly on the noise level, regularization parameter and discretization parameter(s), which gives a useful guideline of the \textsl{a priori} choice of discretization parameters with respect to the noise level in practical implementation. The error analysis is achieved using the conditional stability argument and discrete maximum-norm resolvent estimates. Several numerical experiments are also given to illustrate and complement the theoretical analysis.

\end{abstract}

\section{Introduction}
In this work, we study the inverse problem of recovering a space-dependent diffusion coefficient in 
(sub)diffusion equation from a terminal observation. Let $\Omega\subset \mathbb{R}^d\ (d=1,2,3)$ be a simply connected convex bounded domain with a smooth boundary $\partial\Omega$. The governing equation is given by
\begin{equation}\label{eqn:fde}
   \left\{\begin{aligned}
		\partial_t^\alpha u-\nabla\cdot(q\nabla u) &= f, \ &\mbox{in}&\ \Omega\times(0,T), \\
		u&=0, \ &\mbox{on}&\ \partial\Omega\times(0,T), \\
		u(0)&=u_0, \ &\mbox{in}&\ \Omega,
	\end{aligned}
	\right.
\end{equation}
where $T>0$ is the final time and the notation $\partial_t^\alpha u$ denotes the standard left-sided Djrbashian--Caputo fractional derivative of order $\alpha\in(0,1]$ in the time variable $t$ defined by \cite[p. 92]{kilbas2006theory}:
\begin{equation*}
	\partial_t^\alpha u(t):= \left\{
	\begin{aligned}
		&\frac{1}{\Gamma(1-\alpha)}\int_0^t(t-s)^{-\alpha}\partial_su(s)\ \mathrm{d}s, \ &\mbox{for}&\ \alpha\in(0,1), \\	
		&\partial_tu(t), \ &\mbox{for}&\ \alpha=1,
	\end{aligned}
	\right.
\end{equation*}
with $\Gamma(z)=\int_{0}^{\infty}e^{-s}s^{z-1}\ \mathrm{d}s$, $\Re(z)>0$, being Euler's Gamma function.
The functions $f$ and $u_0$ in \eqref{eqn:fde} are given time-independent source and
initial data, respectively. Due to its extraordinary modeling capability for describing the dynamics of subdiffusion processes
(in which the mean square variance grows sublinearly with the time $t$), the model
\eqref{eqn:fde} has attracted much attention in physics, biology and finance etc.
It has been successfully applied to many important research fields, e.g., subsurface flow
\cite{giona1992fractional,metzler2000random}, thermal diffusion in media with fractal geometry \cite{nigmatullin1986realization}, transport in column experiments \cite{hatano1998dispersive} and highly heterogeneous aquifer \cite{adams1992field}. The classical diffusion model (i.e., $\alpha=1$) represents the most popular to describe transport phenomena found in the nature.

In this work, the concerned inverse problem of the model \eqref{eqn:fde} is to recover the
unknown diffusion coefficient $q^\dag$ from a noisy terminal
observation $z^\delta$:
\begin{equation*}
   z^\delta(x) = u(q^\dag)(x,T) + \xi(x), \quad x\in \Omega,
\end{equation*}
where the exact data $u(q^\dagger)(T)$ denotes the solution of problem \eqref{eqn:fde} (corresponding to $q^\dagger$)
and $\xi$ denotes the pointwise measurement noise. The accuracy of the data $z^\delta$ is measured by the noise level $\delta=\|u(q^\dagger)(T)-z^\delta\|_{L^2(\Omega)}$. The admissible set $\mathcal{A}$ is defined by $ \mathcal{A}=
\{q\in H^1(\Omega): c_0\leq q(x)\leq c_1 \mbox{ a.e. in } \Omega\}$, with $0<c_0<c_1<\infty$. 
Due to the ill-posedness and highly nonlinearity, the numerical recovery of the diffusion coefficient
is challenging.

The study of diffusion coefficient identification in anomalous diffusion has a notable history, dating back to at least the work \cite{cheng2009uniqueness}.  In the one-dimensional case, Cheng et al.  \cite{cheng2009uniqueness} proved the uniqueness for determining a spatially-dependent diffusion coefficient and the fractional order $\alpha$, given the lateral Cauchy data and the Dirac delta function as initial condition. The proof makes use of the Laplace transform and Gel'fand-Leviton theory for inverse Sturm-Liouville problems. Zhang \cite{zhang2015undetermined}
proved the unique recovery of a time dependent diffusion coefficient from lateral Cauchy data. The study with terminal data, despite being highly practical, is still not well understood. Indeed, analyzing terminal data in the standard parabolic case where $\alpha=1$ \cite{Alessandrini:2020} has long been a challenging task, and has only been sparingly studied \cite{AlessandriniVessella:1985,Isakov:1991}. In the one-dimensional case, Isakov \cite{Isakov:1991} analyzed this inverse problem under
some special assumptions on the boundary data. More recently, normal diffusion with a zero source ($f=0$), Triki \cite{Triki:2021} established a Lipschitz stability result
\begin{equation}\label{eqn:Triki}
\| q_1 - q_2  \|_{L^2\II}\le c_T \|   (u(q_1) - u(q_2))(T)\|_{H^2\II}.
\end{equation}
This result holds for sufficiently large $T$ under certain positivity conditions on the initial data $u_0$, achieved by using careful spectral perturbation
estimates (see Remark \ref{rmk:Triki} for further details). However, this spectral perturbation argument in \cite{Triki:2021} is not directly applicable in the error analysis of fully discrete schemes. Moreover, this analysis relies heavily on the exponential decay property of the parabolic problem's solution operator, making it unsuitable for the subdiffusion model with $\alpha\in(0,1)$, where the solution operator decays only linearly. In this paper, we aim to address this gap by proposing a novel conditional stability result for the inverse problem. This result leverages a weighted energy estimate and is applicable for both normal diffusion ($\alpha=1$) and subdiffusion ($\alpha\in(0,1)$).
Additionally, the strategy is amenable with the numerical analysis of discrete schemes.

Numerically, Li et al \cite{li2012numerical,li2013simultaneous} presented the first numerical
recovery of the diffusion coefficient in the fractional case, including both smooth and nonsmooth data, but
without an error analysis of the discrete scheme. Note that in practical computation,
the regularized formulation is often discretized with the Galerkin FEM. The
convergence of discrete approximations as the discretization parameters tend to zero
has been analyzed; See \cite{yamamoto2001simultaneous,keung1998numerical} for the
standard parabolic case. However, deriving a convergence rate is far more challenging,
due to the high degree of nonlinearity of the forward map and strong nonconvexity of
the regularized functional. Thus there have been only very few error bounds on discrete
approximations in the existing literature \cite{wang2010error,JinZhou:2021SINUM,JinZhou:2023IMA,zhang2022identification},
even though such \textit{a priori} estimates can provide useful guidelines for the proper
choice of discretization parameters. The analysis techniques in all these existing works require
that the observational data is available over a time interval (for $\alpha=1$) or whole
space-time interval (for $\alpha\in(0,1)$). The main technical tools include conditional
stability and smoothing properties of solution operator. The current work aims to significantly extend the argument to cover terminal data, which is notably more practical.

In this work, we develop a numerical procedure for recovering the diffusion coefficient
$q$ using a regularized formulation \cite{engl1996regularization,ItoJin:2015} and establish error bounds on
the approximation. We make two new contributions in the work. First, under mild conditions on
the problem data ($u_0$, $f$, $T$, $q_1$, $q_2$ and $\Omega$), we prove a H\"{o}lder type conditional stability in Theorem \ref{thm:main-cont}:
\begin{equation}\label{eqn:stab-intro}
	\|q_1-q_2\|_{L^2(\Omega)}\leq c\|\nabla \big(u(q_1)-u(q_2)\big)(T)\|^\frac12_{L^2(\Omega)}.
\end{equation}
The overall proof relies only on a weighted energy argument (inspired by Bonito et al.
\cite{Bonito:2017}), some nonstandard smoothing properties and asymptotics of solution
operators \cite{Jin:2021book}, and maximum-norm resolvent estimates
\cite{bakaev2003maximum,bakaev2001maximum,stewart1974generation}. 
To the best of our knowledge, this is the first stability result addressing the inverse problem for both the integer-order and fractional-order cases. Moreover, the analysis strategy also plays an essential role in
the error analysis of the inversion scheme.

Second, we employ the standard output least square formulation
to identify the diffusion coefficient.
Motivated by the conditional stability analysis,  an $H^1(\Omega)$-seminorm penalty is used in the formulation. Numerically, both Tikhonov functional and PDE constraint, i.e., problem
\eqref{eqn:fde}, are discretized using the standard Galerkin finite element method
(FEM) with continuous piecewise linear finite elements in space and backward
Euler convolution quadrature in time; see e.g., \cite{lubich1986discretized} and \cite[Chapter 4]{JinZhou:2023book}.
In particular, let $h$ be the spatial mesh size, $\tau$ the time step size, $\gamma$ the regularization parameter, and $q_h^*$ denote the numerical reconstruction of the diffusion coefficient $q^\dag$.
We derive the following error estimate for the numerical approximation in Theorem \ref{thm:main-disc}:
\begin{equation*}
		\|q^\dag-q_h^*\|_{L^2(\Omega)}\leq c(h\gamma^{-1}\eta^2+ \min(1,h^{-1}\eta)\gamma^{-\frac12}\eta+h^2+\tau)^\frac12.
\end{equation*}
with $\eta=\delta + h^2 + \tau + \gamma^\frac12$. The technical proof heavily relies on the conditional stability estimate \eqref{eqn:stab-intro} and
several new smoothing properties and asymptotics of
semi- and fully-discrete solution operators.
Note that the analysis does not involve standard source type conditions,
as is commonly done for nonlinear inverse problem \cite{engl1996regularization,ItoJin:2015}.
The derived $L^2(\Omega)$ error bound is given explicitly in terms of the discretization
parameters $h$ and $\tau$, the noise level $\delta$ and the regularization parameter
$\gamma$ when the fixed value $T$ is relatively large. 
Compared with existing works \cite{wang2010error,JinZhou:2021SINUM,JinZhou:2023IMA,zhang2022identification}, the present work requires overcoming new technical challenges. The key techniques for deriving conditional stability \eqref{eqn:stab-intro} include decay estimate in Lemma \ref{lem:dut} and decay Lipschitz stability in Lemma \ref{lem:dt:u1-u2} of $\partial_t^\alpha u$. Moreover, in the error analysis of the fully discrete scheme, the crucial discrete decay Lipschitz stability estimate does not follow as the continuous case, e.g., maximum-norm resolvent estimates. We develop innovative techniques to overcome the challenge, e.g., the decay estimates of the semi- and fully
discrete solution operators (and their derivatives).
The argument is applicable to both normal diffusion ($\alpha=1$) and subdiffusion ($0<\alpha<1$), thereby significantly broadening the scope of existing works. Numerical experiments indicate that the conditional stability does not hold for small
$T$, cf. Table \ref{tab:err-1d}, confirming the sharpness of the theoretical result.

The rest of the paper is organized as follows. In Section \ref{sec:stab}, we show the conditional stability of the inverse problem. Then in Section \ref{sec:error}, we describe the numerical reconstruction scheme, and provide a complete error analysis for discrete approximations. Finally, in Section \ref{sec:numer}, we present one- and two-dimensional numerical experiments to complement the theoretical results. Last, we give some useful notations. For any $m\geq0$ and $p\geq1$, we denote by $W^{m,p}(\Omega)$ and $W^{m,p}_0(\Omega)$ the standard Sobolev spaces of order $m$, equipped with the norm $\|\cdot\|_{W^{m,p}(\Omega)}$ \cite{Adams2003Sobolev}. We denote by $W^{-m,p'}(\Omega)$ the dual space of $W^{m,p}_0(\Omega)$, with $p'$ being the conjugate exponent of $p$. Further, we write $H^{m}(\Omega)$ and $H^{m}_{0}(\Omega)$ with the norm $\|\cdot\|_{H^m(\Omega)}$ if $p=2$ and write $L^p(\Omega)$ with the norm $\|\cdot\|_{L^p(\Omega)}$ if $m=0$. The notation $(\cdot,\cdot)$ denotes the $L^2(\Omega)$ inner product. We also use Bochner spaces: for a Banach space $B$, let
\begin{equation*}
	W^{m,p}(0,T;B) = \{v: v(\cdot,t)\in B\ \mbox{for a.e.}\ t\in(0,T)\ \mbox{and}\ \|v\|_{W^{m,p}(0,T;B)}<\infty \}.
\end{equation*}
The space $L^\infty(0,T;B)$ is defined similarly. Throughout, we denote by $c$, with or without a subscript, a generic constant which may differ at each occurrence but is always independent of the discretization parameters $h$ and $\tau$, the noise level $\delta$, the regularization parameter $\gamma$ and the terminal time $T$.

\section{Conditional stability} \label{sec:stab}

In this section, we establish a novel conditional stability estimate for the inverse conductivity
problem with the terminal data. Let the operator $A(q)$ be the realization of $-\nabla\cdot(q\nabla\cdot)$ with
a zero Dirichlet boundary condition, with a domain ${\rm Dom}(A(q)):=H^2(\Omega)\cap H_0^1
(\Omega)$. Then the solution $u(q)$ to problem \eqref{eqn:fde} is given by
\begin{equation}\label{eqn:sol}
	u(q) = F(t;q)u_0 + \int_{0}^{t}E(s;q)f\ \mathrm{d}s,
\end{equation}
where the solution operators $F(t;q)$ and $E(t;q)$ are given respectively by \cite[Section 6.2.1]{Jin:2021book}
\begin{equation}\label{eqn:oper-EF}
\begin{split}
	F(t;q)&=\frac{1}{2\pi i}\int_{\Gamma_{\theta,\sigma}}e^{zt}z^{\alpha-1}(z^\alpha+A(q))^{-1}\ \mathrm{d}z
	\quad \mbox{and} \\
	E(t;q)&=\frac{1}{2\pi i}\int_{\Gamma_{\theta,\sigma}}e^{zt}(z^\alpha+A(q))^{-1}\ \mathrm{d}z,
\end{split}
\end{equation}
with the contour $\Gamma_{\theta,\sigma}\subset\mathbb{C}$ (oriented with an increasing imaginary part) given by
$\Gamma_{\theta,\sigma} = \{z\in \mathbb{C}:|z|=\sigma,|\arg(z)|\leq\theta\}\cup\{z\in \mathbb{C}:z=\rho e^{\pm i\theta},\rho\geq\sigma\}$.
Throughout, we fix $\theta\in(\frac{\pi}{2},\pi)$ so that $z^\alpha\in\Sigma_{\alpha\theta}\subset\Sigma_{\theta}$ for all $z\in\Sigma_{\theta}:=\{ z\in\mathbb{C}\setminus\{0\}:|\arg(z)|\leq\theta\}$. The following identity holds:
$\partial_t F(t;q)=-A(q)E(t;q)$ \cite[Lemma 6.2]{Jin:2021book}.

The next lemma gives useful smoothing properties of the operators $F(t;q)$ and $E(t;q)$. For any $s\in\mathbb{R}$, the notation $A(q)^s$ denotes the fractional power of $A(q)$, defined by spectral decomposition.
The cases $s = 0,1$ are known \cite[Theorem 6.4]{Jin:2021book}, and the case $0<s<1$ follows from standard interpolation theory \cite[Proposition 2.3]{LionsMagenes:1972V1}.
\begin{lemma}\label{lem:smoothing properties}
For any $q\in \mathcal{A}$, and any $s\in [0,1]$, there exists $c>0$ independent of $q$ such that
\begin{equation*}
   t^{s\alpha}\|A(q)^sF(t;q)\|_{L^{2}(\Omega)\to L^2(\Omega)} + t^{1-(1-s)\alpha}\|A(q)^sE(t;q)\|_{L^{2}(\Omega)\to L^2(\Omega)}\leq c.
\end{equation*}
\end{lemma}

The following maximum norm resolvent estimate  will be used extensively. See \cite[Theorem 1]{stewart1974generation}, \cite[Theorem 1.1]{bakaev2001maximum} and \cite[Theorem 2.1]{bakaev2003maximum} for the proof.
\begin{lemma}\label{lem:maximum-resolv}
For $q\in W^{1,\infty}(\Omega)\cap \mathcal{A}$, there holds for any $z\in \Sigma_{\theta}$ and $\theta\in(\frac\pi2,\pi)$
\begin{equation}\label{ineq:max-resolvents}
		|z|\|\big(z+A(q)\big)^{-1}\|_{L^\infty(\Omega)\to L^\infty(\Omega)}+|z|^{\frac12}\|\big(z+A(q)\big)^{-1}\|_{L^\infty(\Omega)\to W^{1,\infty}(\Omega)}\leq c.
	\end{equation}
\end{lemma}

Using Lemma \ref{lem:maximum-resolv},
we can derive an \textit{a priori} estimate for the time-fractional derivative $\partial_t^\alpha u$.
\begin{lemma}\label{lem:dut}
Let $u_0\in W^{2,\infty}(\Omega)\cap H_0^1(\Omega)$, $f\in L^\infty(\Omega)$, and $q\in W^{1,\infty}(\Omega)\cap \mathcal{A}$, and let $u(q)$ be the solution to problem \eqref{eqn:fde}. Then there exists $c>0$, independent of $t$ and $q$, such that
\begin{equation*}
   \|\partial_t^\alpha u(q)\|_{W^{1,\infty}(\Omega)}\leq ct^{-\frac{\alpha}{2}}.
\end{equation*}
\end{lemma}
\begin{proof}
By Lemma \ref{lem:maximum-resolv} and letting $\sigma=t^{-1}$ in the contour $\Gamma_{\theta,\sigma}$, we have
\begin{align}\label{ineq:F-inf}
		\|F(t;q)\|_{L^\infty(\Omega)\to W^{1,\infty}(\Omega)}&\leq c\int_{\Gamma_{\theta,\sigma}}|e^{zt}||z|^{\alpha-1}\|\big(z^\alpha+A(q)\big)^{-1}\|_{L^\infty(\Omega)\to W^{1,\infty}(\Omega)}\ |\mathrm{d}z| \nonumber\\
		&\leq c\int_{\Gamma_{\theta,\sigma}}|e^{zt}||z|^{\frac\alpha2-1}\ |\mathrm{d}z|
		\leq ct^{-\frac\alpha2}.
	\end{align}
Let $w(t) = \partial_t^\alpha u(q)\in H_0^1(\Omega)$. Then it satisfies
\begin{equation*}
\left\{\begin{aligned}
	\partial_t^\alpha w-\nabla\cdot(q\nabla w) &= 0, \ &&\mbox{in } \Omega\times(0,T), \\
         w & =0 ,&& \mbox{on } \partial\Omega\times(0,T),\\
			w(0)&=\nabla\cdot(q\nabla u_0)+f, \ &&\mbox{in } \Omega.
		\end{aligned}\right.
	\end{equation*}
It follows from the representation \eqref{eqn:sol} that
$\partial_t^\alpha u(q) = F(t;q)\big(\nabla\cdot(q\nabla u_0)+f\big).$
This, the estimate \eqref{ineq:F-inf} and the assumption on $u_0$ and $f$ lead to
\begin{equation*}
\|\partial_t^\alpha u(q)\|_{W^{1,\infty}(\Omega)}
 \leq c\|F(t;q)\|_{L^\infty(\Omega)\to W^{1,\infty}(\Omega)}\|\nabla\cdot(q\nabla u_0)+f\|_{L^\infty(\Omega)}
\leq ct^{-\frac{\alpha}{2}}.
\end{equation*}
This completes the proof of the lemma.
\end{proof}

Next we provide a crucial Lipschitz stability estimate of the time (fractional) derivative
with respect to the $L^2(\Omega)$ norm of the diffusion coefficient.
\begin{lemma}\label{lem:dt:u1-u2}
Let $u_0\in W^{2,\infty}(\Omega)\cap H_0^1(\Omega)$ and $f\in L^\infty(\Omega)$, and let $u(q_1)$ and $u(q_2)$ be the solutions of problem \eqref{eqn:fde} with $q_1\in \mathcal{A}$ and $q_2\in W^{1,\infty}(\Omega)\cap \mathcal{A}$, respectively. Then for any small $\epsilon>0$,
there exists  $c>0$, independent of $q_1$ and $t$, such that
\begin{equation*}
    \|\partial_t^\alpha(u(q_1)-u(q_2))\|_{L^2(\Omega)}\leq c\max(t^{-\alpha},t^{-\frac{\alpha}{2}(1-2\epsilon)})\|q_1-q_2\|_{L^2(\Omega)},
\end{equation*}
\end{lemma}
\begin{proof}
	Let $u_i = u(q_i)$ and $w(t) = \partial_t^\alpha(u_1-u_2)$. Then $w$  satisfies
	\begin{equation*}
		\left\{
		\begin{aligned}
			\partial_t^\alpha w-\nabla\cdot(q_1\nabla w) &= \nabla\cdot\big((q_1-q_2)\nabla \partial_t^\alpha u_2\big), \ &&\mbox{in}\ \Omega\times(0,T), \\
			w&=0, \ &&\mbox{on}\ \partial\Omega\times(0,T), \\
			w(0)&=\nabla\cdot\big((q_1-q_2)\nabla u_0\big), \ &&\mbox{in}\ \Omega.
		\end{aligned}
		\right.
	\end{equation*}
	The solution representation \eqref{eqn:sol} leads to
	\begin{equation*}
		w(t) = F(t;q_1)\big(\nabla\cdot\big((q_1-q_2)\nabla u_0\big)\big) + \int_{0}^t E(t-s;q_1)\big(\nabla\cdot\big((q_1-q_2)\nabla \partial_s^\alpha u_2(s)\big)\big) \ \mathrm{d}s :={\rm I}_1 + {\rm I}_2.
	\end{equation*}
Since $q_1\in \mathcal{A}$, there exist constants $c$ and $c'$ independent of $q_1$ such that
$$c\|A(q_1)^{\frac12}v\|_{L^2(\Omega)} \le \|\nabla v\|_{L^2(\Omega)}
\le c'\|A(q_1)^{\frac12}v\|_{L^2(\Omega)},\quad \forall v\in H_0^1(\Omega).$$
Then the self-adjointness of the operator $F(t;q_1)$, and integration by parts imply
\begin{align*}
\|{\rm I}_1\|_{L^2(\Omega)} &=\sup_{\|\varphi\|_{L^2(\Omega)}=1}(F(t;q_1)\nabla\cdot\big((q_1-q_2)\nabla u_0,\varphi)\\
  & =\sup_{\|\varphi\|_{L^2(\Omega)}=1}(\nabla\cdot\big((q_1-q_2)\nabla u_0,F(t;q_1)\varphi)\\
  &=\sup_{\|\varphi\|_{L^{2}(\Omega)}=1}\big((q_2-q_1)\nabla u_0,\nabla A(q_1)^{\frac12}F(t;q_1)A(q_1)^{-\frac12}\varphi\big).
\end{align*}
Then Lemma \ref{lem:smoothing properties} and the boundedness of the operator
$A(q_1)^{-\frac12}$ in $L^2(\Omega)$ (uniform in $q_1$) imply
\begin{align*}
\|{\rm I}_1\|_{L^2(\Omega)} \leq &c\|q_1-q_2\|_{L^2(\Omega)}\|\nabla u_0\|_{L^\infty(\Omega)}
\| A(q_1)F(t;q_1)\|_{L^2(\Omega) \to L^2(\Omega)}\sup_{\|\varphi\|_{L^{2}(\Omega)}=1}\|A(q_1)^{-\frac12}\varphi\|_{L^2(\Omega)}\\
\leq & ct^{-\alpha}\|q_1-q_2\|_{L^2(\Omega)}.
\end{align*}
Similarly, by the boundedness of the operator $A(q_1)^{-\frac12+\epsilon}$ in $L^2(\Omega)$ (uniform in $q_1$) for small $\epsilon>0$ and Lemmas \ref{lem:smoothing properties}-\ref{lem:dut}, we have
\begin{align*}
&\|E(t-s;q_1)\big(\nabla\cdot\big((q_1-q_2)\nabla \partial_s^\alpha u_2(s)\big)\big)\|_{L^2(\Omega)} \\ =&\sup_{\|\varphi\|_{L^{2}(\Omega)}=1}\big((q_2-q_1)\nabla \partial_s^\alpha u_2(s),\nabla A(q_1)^{\frac12-\epsilon}E(t-s;q_1)A(q_1)^{-\frac12+\epsilon}\varphi\big)\\
\leq& c\|q_1-q_2\|_{L^2(\Omega)}\|\nabla \partial_s^\alpha u_2(s)\|_{L^\infty(\Omega)}
\| A(q_1)^{1-\epsilon}E(t-s;q_1)\|_{L^2(\Omega) \to L^2(\Omega)}\sup_{\|\varphi\|_{L^{2}(\Omega)}=1}\|A(q_1)^{-\frac12+\epsilon}\varphi\|_{L^2(\Omega)}\\
\leq &c(t-s)^{\epsilon\alpha-1}s^{-\frac{\alpha}{2}}\|q_1-q_2\|_{L^2(\Omega)}.
\end{align*}
Consequently,
\begin{align*}
\|{\rm I}_2\|_{L^2(\Omega)} & \leq  \int_{0}^t\|E(t-s;q_1)\big(\nabla\cdot((q_1-q_2)\nabla \partial_s^\alpha u_2(s))\big)\|_{L^2(\Omega)} \ \mathrm{d}s\\
& \leq c\|q_1-q_2\|_{L^2(\Omega)} \int_{0}^t(t-s)^{\epsilon\alpha-1}s^{-\frac{\alpha}{2}}\ \mathrm{d}s\leq ct^{-\frac{\alpha}{2}(1-2\epsilon)}\|q_1-q_2\|_{L^2(\Omega)}.
\end{align*}
The bounds on ${\rm I}_1$ and ${\rm I}_2$ and the triangle inequality complete the proof of the lemma.
\end{proof}

Next we give a novel conditional stability estimate. First we state the standing assumption.
\begin{assumption}\label{ass:cont}
$q^\dag\in W^{1,\infty}(\Omega)\cap \mathcal{A}$, $u_0\in W^{2,\infty}(\Omega)\cap H_0^1(\Omega)$, and $f\in L^{\infty}(\Omega)$.
\end{assumption}

Under Assumption \ref{ass:cont}, the following regularity results hold.
Let $p>\max(d,2)$. Then for any $\theta<\frac12-\frac{d}{2p}$ and $r>\frac{1}{\alpha\theta}$,
the solution $u = u(q^\dag)$ to problem \eqref{eqn:fde} satisfies \cite[(2.5)--(2.6)]{JinZhou:2023IMA}
\begin{enumerate}
 	\item[(i)] $u\in W^{\alpha\theta,r}(0,T;W^{2(1-\theta),p})\hookrightarrow L^{\infty}(0,T;W^{1,\infty}(\Omega))$;
 	\item[(ii)] $\|u(t)\|_{H^2(\Omega)}+\|\partial_t^\alpha u(t)\|_{L^2(\Omega)}+t^{1-\alpha}\|\partial_t u(t)\|_{L^2(\Omega)}+t\|\partial_tu(t)\|_{L^2(\Omega)}\leq c$, a.e. $t\in(0,T]$.
\end{enumerate}

Additionally, we assume that the following positivity condition holds:
\begin{equation}\label{Cond:P}
(q^\dag|\nabla u(q^\dag)|^2+(f-\partial_t^\alpha u(q^\dag))u(q^\dag))(T) \ge c >0.
\end{equation}
This condition was proved in \cite{JinZhou:2021SINUM} for parabolic equations (i.e., $\alpha=1$)
and in \cite{JinZhou:2023IMA} for time-fractional diffusion (i.e., $0<\alpha<1$). For example, if $\Omega$ is a
$C^{2,\mu}$ domain with $\mu\in(0,1)$, $q^\dag\in C^{1,\mu}(\Omega)\cap \mathcal{A}$,
$f\in C^\mu(\Omega)$ with $f\geq c_f>0$ and $u_0\in C^{2,\mu}(\Omega)\cap H_0^1(\Omega)$ with
$u_0\geq0$ in $\Omega$, and $f+\nabla\cdot(q^\dag\nabla u_0)\leq0$ in $\Omega$, then condition \eqref{Cond:P}
holds. See \cite[Section 4.3]{JinZhou:2021SINUM} and \cite[Proposition 3.5]{JinZhou:2023IMA} for detailed discussions.

Now we give a H\"{o}lder type conditional stability estimate for the inverse problem. The estimate is
conditional since the coefficients are required to have extra regularity. To the best of our knowledge, this appears to be the first result of the kind for
the time-fractional model \eqref{eqn:fde} with a terminal observation. The analysis strategy will also guide the error
analysis of the fully discrete scheme in Section \ref{sec:error} below.  
\begin{theorem}\label{thm:main-cont}
Let Assumption \ref{ass:cont} hold, and $q\in \mathcal{A}$ with $\| \nabla q \|_{L^2(\Omega)} \le c$.
Then for small $\epsilon>0$, the following conditional stability estimate holds
\begin{align*}
	&\int_{\Omega}\Big(\frac{q^\dag-q}{q^\dag}\Big)^2(q^\dag|\nabla u(q^\dag)|^2+(f-\partial_t^\alpha u(q^\dag))u(q^\dag))(T)\ \mathrm{d}x\\
 \leq& c\|\nabla\big(u(q) -u (q^\dag)\big)(T)\|_{L^2(\Omega)} + c\max(T^{-\alpha},T^{-\frac{\alpha}{2}(1-2\epsilon)})\|q^\dag-q\|^2_{L^2(\Omega)},
\end{align*}
with $c>0$ independent of $q$ and $T$. Moreover, if condition \eqref{Cond:P} holds,
then there exist $T_0>0$ and $c>0$, independent of $q$, such that for all $T\ge T_0$,
\begin{equation*}
	\|q-q^\dag\|_{L^2(\Omega)}\leq c\|\nabla\big(u(q) -u(q^\dag) \big)(T)\|^{\frac12}_{L^2(\Omega)}.
\end{equation*}
\end{theorem}
\begin{proof}
Let $u^\dag=u(q^\dag)$ and $ u=u(q)$. Then by the weak formulations of $u^\dag$ and
$u$, there holds for any $\varphi\in H_0^1(\Omega)$
\begin{equation*}
	\big((q^\dag- q)\nabla u^\dag(T),\nabla\varphi\big) = -\big(q\nabla(u^\dag- u)(T),\nabla\varphi\big)-\big(\partial_t^\alpha(u^\dag- u)(T),\varphi\big)=: {\rm I_1}+{\rm I}_2.
\end{equation*}
Let $\varphi\equiv\frac{q^\dag- q}{q^\dag}u^\dag(T)\in H^1_0(\Omega)$. Upon repeating the argument
in \cite{Bonito:2017,JinZhou:2021SINUM}, we obtain
    \begin{equation*}
    	\big((q^\dag- q)\nabla u^\dag(T),\nabla\varphi\big) = \frac12\int_{\Omega}\Big(\frac{q^\dag- q}{q^\dag}\Big)^2\big(q^\dag|\nabla u^\dag|^2+(f-\partial_t^\alpha u^\dag)u^\dag\big)(T)\ \mathrm{d}x.
    \end{equation*}
Meanwhile, direct computation yields
     \begin{equation*}
     	\nabla\varphi
     	=\nabla\Big(\frac{q^\dag- q}{q^\dag}\Big) u^\dag(T)+\Big(\frac{q^\dag- q}{q^\dag}\Big)\nabla u^\dag(T).
     \end{equation*}
Now Assumption \ref{ass:cont} implies $\|u^\dag(T)\|_{L^\infty(\Omega)}+\|\nabla u^\dag(T)\|_{L^2(\Omega)}\leq c$. This,
the box constraint on $q^\dag,q\in \mathcal{A}$ and the \textit{a priori} bound $\| \nabla q \|_{L^2(\Omega)} \le c$ yield
\begin{equation*}
   \|\varphi\|_{L^2(\Omega)}\leq c\|q-q^\dag\|_{L^2(\Omega)}\quad \mbox{and}\quad \|\nabla \varphi\|_{L^2(\Omega)} \leq c.
\end{equation*}
Hence, by Lemma \ref{lem:dt:u1-u2} and the Cauchy-Schwarz inequality, we obtain
     \begin{align*}
     	|{\rm I}_1|&\leq c\|\nabla (u-u^\dag)(T)\|_{L^2(\Omega},\\
     	|{\rm I}_2|
      &\leq c\|\partial_t^\alpha (u-u^\dag)(T)\|_{L^2(\Omega)}\|q-q^\dag\|_{L^2(\Omega)}\\
      &\leq c\max(T^{-\alpha},T^{-\frac{\alpha}{2}(1-2\epsilon)})\|q-q^\dag\|^2_{L^2(\Omega)}.
     \end{align*}
Moreover, under the condition \eqref{Cond:P}, we apply the box constraint on $q,q^\dag\in\mathcal{A}$ and derive
\begin{equation*}
     \|q-q^\dag\|^2_{L^2(\Omega)}\leq c\|\nabla\big(u-u^\dag\big)(T)\|_{L^2(\Omega)} + c\max(T^{-\alpha},T^{-\frac{\alpha}{2}(1-2\epsilon)})\|q-q^\dag\|^2_{L^2(\Omega)}.
\end{equation*}
Let $T_0$ be sufficiently large such that
$c\max(T^{-\alpha},T^{-\frac{\alpha}{2}(1-2\epsilon)}) \le \frac12$. Then for any $T\ge T_0$, the desired estimate follows.
\end{proof}

\begin{remark}\label{rmk:Triki}
Theorem \ref{thm:main-cont} extends several existing works. The weighted stability results of the type
for observation over a space-time domain were implicitly obtained in \cite{JinZhou:2021SINUM,
JinZhou:2023IMA}. The only result for
the terminal data case was obtained by Triki \cite[Theorem 1.1]{Triki:2021} for the standard parabolic problem, who proved the following Lipschitz stability for a large $T$: for $q,q^\dag\in C^1(\overline{\Omega})$, there holds
\begin{equation*}
  \|q-q^\dag\|_{L^2(\Omega)}\leq c\|u(q^\dag)(T)-u(q)(T)\|_{H^2(\Omega)},
\end{equation*}
where the constant $c$ depends on the terminal time $T$ (exponentially) and the domain $\Omega$.
This result was shown for the case $f\equiv0$ and $u_0$ satisfying a mild positivity condition,
and the proof relies on the decay estimate on $\partial_t u$, which itself was
proved using refined spectral perturbation estimates. Theorem \ref{thm:main-cont} provides a novel H\"older stability estimate using an energy estimate and can be adapted to the error analysis of numerical approximations in Section \ref{sec:error}.
\end{remark}

\section{Numerical approximation and error analysis} \label{sec:error}

Now we develop a numerical procedure based on the regularized output least-squares formulation, and discretize the regularized problem using backward Euler convolution quadrature (BECQ) in time and Galerkin FEM with continuous piecewise linear elements in space. Furthermore, we provide a complete error analysis of the fully discrete scheme.

\subsection{Regularized problem and numerical approximation}
\label{subsec:Tikhonov regularization problem and its FEM approximation in elliptic system}
To identify the diffusion coefficient $q$, we employ the standard Tikhonov regularization with an $H^1(\Omega)$ seminorm penalty \cite{engl1996regularization,ItoJin:2015}, which gives the following minimization problem:
\begin{equation}\label{eqn:obj}
	\min_{q\in \mathcal{A}} J_{\gamma}(q)=\frac{1}{2}\|u(q)(T)-z^{\delta}\|^2_{L^2(\Omega)}+\frac{\gamma}{2}\|\nabla q\|_{L^2(\Omega)}^2,
\end{equation}
where $\gamma>0$ is the regularization parameter and $u(t)\equiv u(q)(t)\in H^1_0(\Omega)$ with $u(0)=u_0$ satisfies
\begin{equation}\label{eqn:weak}	
	(\partial^\alpha_tu(t),\varphi)+(q\nabla u(t),\nabla\varphi)=(f,\varphi), \quad\forall \varphi\in H^1_0(\Omega),\ \mbox{a.e.}\ t\in(0,T).
\end{equation}
By a standard argument \cite{engl1996regularization,ItoJin:2015}, it can be proved
that problem \eqref{eqn:obj}--\eqref{eqn:weak} has at least one global minimizer $q_\gamma^\delta$, which is
continuous with respect to the perturbations in the data $z^\delta$. Moreover, as the noise level
$\delta\to0^+$, the sequence $\{q_\gamma^\delta\}_{\delta>0}$ of minimizers contains a subsequence that converges to the exact coefficient $q^\dag$ in $H^1(\Omega)$ if 
$\gamma$ is chosen properly.

In practice, one needs to discretize the regularized formulation \eqref{eqn:obj}--\eqref{eqn:weak}  suitably.
For time discretization, we divide the interval $[0,T]$ uniformly into
$N$ subintervals, with a time step size $\tau:= T/N$ and grid $t_n:=n\tau$,
$n=0,1\dots N$. To approximate the fractional derivative $\partial_t^\alpha
v(t_n)$, we employ backward Euler convolution quadrature (BECQ) defined by
\begin{equation*}
	\bar{\partial}_\tau^\alpha v^n :=\tau^{-\alpha}\sum_{j=0}^nb^{(\alpha)}_{j}(v^{n-j}-v^0), \quad \mbox{with}  \ v^j=v(t_j),
\end{equation*}
where the weights $b_j^{(\alpha)}$ are generated by the power series expansion $(1-\zeta)^\alpha=\sum_{j=0}^\infty b_j^{(\alpha)}\zeta^j$. The weights $b_j^{(\alpha)}$ are given by $b_j^{(\alpha)}=(-1)^j\alpha(\alpha-1)\cdots(\alpha-j+1)/j!$,
with $b_0^{(\alpha)}=1$ and $b_j^{(\alpha)}<0$ for $j\geq1$. When $\alpha=1$, it reduces to the standard backward Euler scheme.

For the spatial discretization, we employ the standard Galerkin FEM. Let $h\in(0,h_0]$ for some $h_0>0$ and $\mathcal{T}_h:=\cup\{T_j\}_{j=1}^{N_h}$ be a shape regular quasi-uniform simplicial triangulation of the domain $\Omega$ into mutually disjoint open {face-to face 
subdomains $T_j$, such that $\Omega_h:= {\rm Int}(\cup_j\{\bar{T}_j\})\subset\Omega$ with all the
boundary vertices of the domain $\Omega_h$ locating on $\partial\Omega$ and ${\rm dist}(x,\partial\Omega)\leq ch^2$ for $x\in\partial\Omega_h$ \cite[Section 5.3]{LarssonThomee:2003}. On the triangulation $\mathcal{T}_h$, we define the space $V_h$  of continuous piecewise linear finite element functions by
\begin{equation*}
	V_{h}:=
	\{v_{h}\in H^1(\Omega_h):	v_{h}|_{T}\text{ is a linear polynomial},\, \forall T\in\mathcal{T}_h \}.
\end{equation*}
Note that the functions in $V_{h}$ can be naturally extended to the entire domain $\Omega$ by linear polynomials, and we denote the space of extended functions also by $V_h$. Moreover, we define
the space $X_h$  (that vanish outside $\Omega_h$) by
\begin{equation*}
X_{h}:=	\{v_{h}\in H_0^1(\Omega_h):	v_{h}|_{T}\text{ is a linear polynomial }\forall T\in\mathcal{T}_h~\text{and}~v_h|_{\Omega\backslash\Omega_h} = 0 \}
\end{equation*}
The spaces $V_h$ and $X_h$  are used to discretize the diffusion coefficient $q$ and the state $u$, respectively. Note that if $\Omega$ is a convex polygon, then $X_h = V_h\cap H_0^1(\Omega)$. Next
we recall several useful estimates.
We denote by $\Pi_h$ the Lagrange nodal interpolation operator on $V_h$. Since ${\rm dist}(x,\partial\Omega)\leq ch^2$ for $x\in\partial\Omega_h$, we have
(see \cite[Theorem 4.4.20]{BrennerScott:book2008} for a
convex polyhedral domain and Lemma \ref{app:lem} for a convex domain with a curved boundary):
\begin{align}
\|v-\Pi_hv\|_{L^2(\Omega)} + h\|\nabla(v-\Pi_hv)\|_{L^2(\Omega)}  &\leq ch^2 \|v\|_{H^2(\Omega)}, \quad \forall v\in H^2(\Omega),
\label{inequ:Pi_h approx-2}\\
\|v-\Pi_hv\|_{L^\infty(\Omega)} + h\|\nabla(v-\Pi_hv)\|_{L^\infty(\Omega)}  &\leq ch \|v\|_{W^{1,\infty}(\Omega)}, \quad \forall v\in W^{1,\infty}(\Omega).\label{inequ:Pi_h approx-inf}
\end{align}
Moreover, we define the standard $L^2(\Omega)$-projection operator $P_h:L^2(\Omega)\to X_h$ by
\begin{equation*}
	(P_hv,\varphi_h)=(v,\varphi_h), \quad \forall v\in L^2(\Omega),  \varphi_h\in X_h,
\end{equation*}
Then  for $1\leq p\leq \infty$ and $s=0,1,2$, $k=0,1$ with $k\leq s$ \cite{CrouzeixThomee:1987, bakaev2001maximum}:
\begin{equation}\label{inequ:P_h approx}
	\|v-P_hv\|_{W^{k,p}(\Omega)}\leq  Ch^{s-k}\|v\|_{W^{s,p}(\Omega)}, \quad \forall v\in W^{s,p}(\Omega) \cap H_0^1(\Omega).
\end{equation}

Now we can formulate a fully discrete scheme for the regularized problem \eqref{eqn:obj}-\eqref{eqn:weak} as
\begin{equation}\label{eqn:obj-disc}
	\min_{q_h\in \mathcal{A}_h} J_{\gamma,h,\tau}(q_h)=\frac12\|U^N_h(q_h)-z^{\delta}\|^2_{L^2(\Omega)}+\frac{\gamma}{2}\|\nabla q_h\|_{L^2(\Omega)}^2,
\end{equation}
with $\mathcal{A}_h = \mathcal{A}\cap V_h$, where $U_h^n \equiv U^n_h(q_h)\in X_h$ satisfies $U_h^0=P_hu_0$ and
\begin{equation}\label{eqn:fem}	
	(\bar{\partial}_\tau^\alpha U_h^n,\varphi_h)+(q_h\nabla U^n_h,\nabla\varphi_h)=(f,\varphi_h), \quad\forall \varphi_h\in X_h,\ n=1,2,\ldots,N.
\end{equation}
The discrete problem \eqref{eqn:obj-disc}--\eqref{eqn:fem} is well-posed: there exists at
least one global minimizer $q_h^*\in \mathcal{A}_h$, and it depends continuously on the
data. Further, as the discretization parameters $h$ and $\tau$ tend to zero,
the numerical approximation $q_h^*$ converges to the regularized solution to problem
\eqref{eqn:obj}--\eqref{eqn:weak}. We aim to establish a bound on the error $q_h^*-q^\dag$
in terms of the noise level $\delta$, discretization parameters $h$ and $\tau$ and regularization parameter $\gamma$.
For the error analysis, we need the following assumption on the problem data.
\begin{assumption}\label{ass:disc}
$q^\dag\in H^2(\Omega)\cap W^{1,\infty}(\Omega)\cap \mathcal{A}$, $u_0\in W^{2,\infty}(\Omega)\cap H_0^1(\Omega)$, and $f\in L^{\infty}(\Omega)$.
\end{assumption}

Now we give the main result in this section, i.e., a weighted $L^2(\Omega)$
error bound on the approximation $q_h^*$.
The proof heavily relies on some technical estimates, whose proofs are deferred to Section \ref{ssec:est}.
\begin{theorem}\label{thm:main-disc}
Let Assumption \ref{ass:disc} hold, and $\{(q_h^*,u_h^n(q_h^*))\}_{n=0}^N$ be the solutions of problem \eqref{eqn:obj-disc}-\eqref{eqn:fem}. Then with $\eta_T=T^{\alpha-1}\tau+\max(1,T^{-\alpha})h^2+\delta+\gamma^\frac12$ there holds  
\begin{align*}
    &\int_{\Omega}\Big(\frac{q^\dag-q_h^*}{q^\dag}\Big)^2(q^\dag|\nabla u(q^\dag)|^2+(f-\partial_t^\alpha u(q^\dag))u(q^\dag))(T)\ \mathrm{d}x\\ \leq&c\big(h\gamma^{-1}\eta_T^2+ \min(1,h^{-1}\eta_T)\gamma^{-\frac12}\eta_T +  h^2\max(T^{-\alpha}, T^{-2\alpha}) + \tau T^{-\alpha-1}\big) \\
    & \qquad \qquad +c\max(T^{-\alpha}, T^{-\frac\alpha2(1-2\epsilon)},T^{-\alpha(1-\epsilon)},T^{-\alpha(2-\epsilon)}) \|q^\dag-q^*_h\|^2_{L^2(\Omega)}.
\end{align*}
where the generic constants  are independent of $h$, $\tau$, $\delta$, $\gamma$ and $T$. Moreover, under condition \eqref{Cond:P}, with $\eta:=\tau+h^2+\delta+\gamma^\frac12$, there exits $T_0>0$ such that for any $T\ge T_0$,
	\begin{equation*}
		\|q^\dag-q_h^*\|^2_{L^2(\Omega)}\leq c(h\gamma^{-1}\eta^2+ \min(1,h^{-1}\eta)\gamma^{-\frac12}\eta+h^2+\tau).
	\end{equation*}
\end{theorem}
\begin{proof}
The proof proceeds similarly to the conditional stability estimate in Theorem \ref{thm:main-cont}
and requires several (new) technical estimates proved in the propositions below.
Let $u^\dag\equiv u(q^\dag)$. For any test function $\varphi\in H_0^1(\Omega)$, using
the weak formulations of $u^\dag$ and $U_h^N(q_h^*)$, we have
\begin{align*}
&\quad\,  \big((q^\dag-q_h^*)\nabla u^\dag(T),\nabla\varphi\big)\\
   & = \big ((q^\dag-q_h^*)\nabla u^\dag(T),\nabla(\varphi-P_h\varphi)\big) + \big( (q^\dag-q_h^*)\nabla u^\dag(T),\nabla P_h\varphi\big) \\
		& = \big((q^\dag-q_h^*)\nabla u^\dag(T),\nabla(\varphi-P_h\varphi)\big)+(q^*_h\nabla(U_h^N(q_h^*)-u^\dag(T)),\nabla P_h\varphi\big)\\
  &\quad+\big(q^\dag \nabla u^\dag(T)-q_h^*\nabla U_h^N(q_h^*),\nabla P_h\varphi\big) \\
		& = -\big(\nabla\cdot((q^\dag-q_h^*)\nabla u^\dag(T)),\varphi-P_h\varphi\big)+(q^*_h\nabla(U_h^N(q_h^*)-u^\dag(T)),\nabla P_h\varphi\big)\\
       &\quad+\big(\bar{\partial}_\tau^\alpha U_h^N(q_h^*)-\partial_t^\alpha u^\dag(T), P_h\varphi\big)=: {\rm I}_1+{\rm I}_2+{\rm I}_3.
	\end{align*}
Now we bound the three terms ${\rm I}_i$, $i=1,2,3$, separately. Let $\varphi=\frac{q^\dag-q_h^*}{q^\dag}
u^\dag(T)$. Then by the box constraint $q^\dag,q_h^*\in\mathcal{A}$ and Proposition
\ref{prop:err:uT-uhN-and-gqh} below, $\varphi$ satisfies
\begin{equation}\label{eqn:est-varphi-h}
  \|\varphi\|_{L^2(\Omega)}\leq c, \quad
  \|P_h\varphi\|_{L^2(\Omega)} 
  \leq c \quad \mbox{and}\quad \|\nabla\varphi\|_{L^2(\Omega)}\leq c(1+\|\nabla q_h^*\|_{L^2(\Omega)}).
\end{equation}
By Assumption \ref{ass:disc}, we have
\begin{equation*}
  \|\Delta u^\dag(T)\|_{L^2(\Omega)}+\|\nabla u^\dag\|_{L^\infty(\Omega)}\leq c.
\end{equation*}
Thus, direct computation yields
\begin{equation*}
\begin{split}
   &\quad\,\, \|\nabla\cdot\big((q^\dag-q_h^*)\nabla u^\dag(T)\big)\|_{L^2(\Omega)}\\
   &\leq\|q^\dag-q_h^*\|_{L^\infty(\Omega)}\|\Delta u^\dag(T)\|_{L^2(\Omega)}+\|\nabla q^\dag\|_{L^\infty(\Omega)}\|\nabla u^\dag(T)\|_{L^2(\Omega)}\\
   &\quad +\|\nabla q^*_h\|_{L^2(\Omega)}\|\nabla u^\dag(T)\|_{L^\infty(\Omega)}\leq c(1+\|\nabla q_h^*\|_{L^2(\Omega)}).
\end{split}
\end{equation*}
This, the estimate \eqref{eqn:est-varphi-h} and Proposition \ref{prop:err:uT-uhN-and-gqh} below imply
\begin{equation*}
   |{\rm I_1}|\leq ch(1+\|\nabla q_h^*\|_{L^2(\Omega)})\|\nabla \varphi\|_{L^2(\Omega)}\leq ch(1+\|\nabla q_h^*\|^2_{L^2(\Omega)})\leq ch\gamma^{-1}\eta_T^2.
\end{equation*}
Next, by the triangle inequality, the inverse inequality in the space $X_h$ \cite[equation (1.12)]{Thome2006GalerkinFE},
Proposition \ref{prop:err:uT-uhN-and-gqh} below and the approximation property of of $P_h$ in \eqref{inequ:P_h approx}, we have
\begin{align*}
|{\rm I_2}|&\leq \big(\|\nabla(U_h^N(q_h^*)-P_hu^\dag(T))\|_{L^2(\Omega)}+\|\nabla(u^\dag(T)-P_hu^\dag(T))\|_{L^2(\Omega)}\big)\|\nabla\varphi\|_{L^2(\Omega)}\\
&\leq c\big(h^{-1}\| U_h^N(q_h^*)-P_hu^\dag(T)\|_{L^2(\Omega)}+\|\nabla(u^\dag(T)-P_hu^\dag(T))\|_{L^2(\Omega)}\big)\|\nabla\varphi\|_{L^2(\Omega)}\\
&\leq c(h+h^{-1}\eta_T)\gamma^{-\frac12}\eta_T.
\end{align*}
This and the \textit{a priori} estimate $\|\nabla\big(U_h^N(q_h^*)-u^\dag(T)\big)\|_{L^2(\Omega)}\leq c$ imply
\begin{equation*}
	|{\rm I_2}|\leq c\min(1,h^{-1}\eta_T)\gamma^{-\frac12}\eta_T.
\end{equation*}
Next, to bound the term ${\rm I}_3$, we employ the splitting
\begin{align*}
&\big(\bar{\partial}_\tau^\alpha U_h^N(q_h^*)-\partial_t^\alpha u^\dag(T), P_h\varphi\big)\\
=&\big(\bar{\partial}_\tau^\alpha U_h^N(q_h^*)-\bar{\partial}_\tau^\alpha U_h^N(q^\dag), P_h\varphi\big)+\big(\bar{\partial}_\tau^\alpha U_h^N(q^\dag)-\bar{\partial}_\tau^\alpha u^N(q^\dag), P_h\varphi\big)\\&+\big(\bar{\partial}_\tau^\alpha u^N(q^\dag)-\partial_t^\alpha u^\dag(T), P_h\varphi\big)=:{\rm I}_{3}^1+{\rm I}_{3}^2+{\rm I}_{3}^3,
\end{align*}
and then bound the three terms separately. It follows from Lemmas
\ref{lem: err on frac discretization}-\ref{lem: err on space discretization}
below and the estimate \eqref{eqn:est-varphi-h} that
\begin{equation*}
 |{\rm I}_{3}^2|\leq ch^2\max(T^{-\alpha}, T^{-2\alpha})\quad\mbox{and}\quad|{\rm I}_3^3|\leq c\tau T^{-\alpha-1}.
\end{equation*}
Meanwhile, Proposition \ref{prop:err-dtauuq-uqh} below and the estimate \eqref{eqn:est-varphi-h} lead to
\begin{equation*}
|{\rm I}_{3}^1|\leq c\max(T^{-\alpha},T^{-\frac\alpha2(1-2\epsilon)},T^{-\alpha(1-\epsilon)}, T^{-\alpha(2-\epsilon)})\|q^\dag-q^*_h\|^2_{L^2(\Omega)}.
\end{equation*}
Upon combing these estimates with the identity \cite{Bonito:2017,JinZhou:2021SINUM}
\begin{equation*}
	\big((q^\dag-q_h^*)\nabla u^\dag(T),\nabla\varphi\big)=\frac12\int_{\Omega}\Big(\frac{q^\dag-q_h^*}{q^\dag}\Big)^2(q^\dag|\nabla u^\dag|^2+(f-\partial_t^\alpha u^\dag)u^\dag)(T)\ \mathrm{d}x,
\end{equation*}
we prove the first assertion. Since $\eta_T\leq c\eta$ for large $T$, the second assertion follows exactly as Theorem \ref{thm:main-cont}.
\end{proof}

\begin{remark}
The estimate in Theorem \ref{thm:main-disc} provides useful guidelines for choosing the algorithmic parameters:
Given the noise level $\delta$, we may choose $\gamma  \sim  \delta^2$ and $h \sim  \delta^\frac12$. The choice $\gamma  \sim  \delta^2$ differs from the usual condition for Tikhonov
regularization, i.e., $\lim_{\delta\rightarrow 0^+}\frac{\delta^2}{\gamma (\delta )} = 0$,
but it agrees with that with conditional stability (see, e.g., \cite[Theorems 1.1 and
1.2]{EggerHofmann:2018}). It is noteworthy that the error bound in Theorem \ref{thm:main-disc} is comparable with that for the standard parabolic
case \cite[Theorem 4.5]{JinZhou:2021SINUM} and the time fractional case \cite{JinZhou:2021sicon,JinZhou:2023IMA}.
Theorem \ref{thm:main-disc} requires only the terminal data, whereas previous results
\cite{JinZhou:2021sicon,JinZhou:2023IMA} in the fractional case require full space-time data. Thus it represents a substantial improvement for the concerned inverse problem.
\end{remark}

\subsection{Preliminary technical estimates}\label{ssec:est}
In this part, we derive crucial \textit{a priori} bounds on $u(q^\dag)(T)-U_h^N(q^*_h)$ and $\nabla q_h^*$;
see Proposition \ref{prop:err:uT-uhN-and-gqh} for the precise statement.
For any $q\in \mathcal{A}$, we define a discrete elliptic operator
$A_h(q): X_h\to X_h$ by
\begin{equation*}
	(A_h(q)v_h,\varphi_h)=(q\nabla v_h, \nabla \varphi_h),\quad \forall v_h,\varphi_h\in X_h.
\end{equation*}
Then problem \eqref{eqn:fem} is equivalent to an operator equation in $X_h$:
\begin{equation*}
	\bar{\partial}_\tau^\alpha U_h^n + A_h(q_h)U_h^n = P_hf, \quad\mbox{with}\ U_h^0=P_hu_0,\ n=1,2,\dots,N.
\end{equation*}
Using the discrete Laplace transform, the solution $U_h^n(q_h)$ is given by
\begin{equation*}
	U_h^n(q_h)=F_{h,\tau}^n(q_h)U_h^0 + \tau\sum_{j=1}^nE_{h,\tau}^j(q_h)P_hf=F_{h,\tau}^n(q_h)P_hu_0 + (I-F_{h,\tau}^n(q_h))A_h(q_h)^{-1}P_hf,
\end{equation*}
where the fully discrete solution operators $F_{h,\tau}^n(q_h)$ and $E_{h,\tau}^n(q_h)$ are defined respectively by \cite[Section 3.2]{JinZhou:2023book}
\begin{equation}\label{equ: dis oper F and E}
\begin{aligned}
	F_{h,\tau}^n(q_h)&=\frac{1}{2\pi i}\int_{\Gamma^\tau_{\theta,\sigma}}e^{zt_{n-1}}\delta_\tau(e^{-z\tau})^{\alpha-1}(\delta_\tau(e^{-z\tau})^{\alpha}+A_h(q_h))^{-1}\ \mathrm{d}z,  \\
	E_{h,\tau}^n(q_h)&=\frac{1}{2\pi i}\int_{\Gamma^\tau_{\theta,\sigma}}e^{zt_{n-1}}(\delta_\tau(e^{-z\tau})^{\alpha}+A_h(q_h))^{-1}\ \mathrm{d}z,
\end{aligned}
\end{equation}
with the kernel function $\delta_\tau(\zeta):=\tau^{-1}(1-\zeta)$ and the contour
$\Gamma_{\theta,\sigma}^\tau:=\{z\in\Gamma_{\theta,\sigma}:|\Im(z)|\leq \frac{\pi}{\tau}\}$ (oriented with an increasing imaginary part).
Like in the continuous case, we need suitable smoothing properties of the fully discrete solution operators $F_{h,\tau}^n(q)$ and $E_{h,\tau}^n(q)$. Recall that for any fixed $\theta\in (\frac{\pi}{2},\pi)$, there exists $\theta'\in(\frac{\pi}{2},\pi)$ such that for all $\alpha\in(0,1]$ and $z\in \Gamma^\tau_{\theta,\sigma}$ \cite[Lemma 3.1]{JinZhou:2023book}:
\begin{equation}\label{ineq: dis kernel property}
	c_1|z|\leq|\delta_\tau(e^{-z\tau})|\leq c_2|z|,\quad\delta_\tau(e^{-z\tau})\in\Sigma_{\theta'}, \quad|\delta_\tau(e^{-z\tau})^\alpha-z^\alpha|\leq c_3\tau|z|^{1+\alpha},
\end{equation}
where the constants $c_1$, $c_2$ and $c_3$ are independent of $\tau$. Further, for any
fixed $q\in \mathcal{A}$, let $\lambda(q)$ and $\lambda_h(q)$ be the smallest
eigenvalues of operators $A(q)$ and $A_h(q)$, respectively. Then by Courant-Fischer-Weyl minmax
theorem, there exists $c_0>0$, independent of $q$, such that $c_0\leq\lambda(q)
\leq\lambda_h(q)$. The following discrete resolvent estimate holds
\begin{equation}\label{ineq: dis resol}
\begin{split}
	\|\big(\delta_\tau(e^{-z\tau})^\alpha+A_h(q)\big)^{-1}\|_{L^2(\Omega)\to L^2(\Omega)} & \leq c\min(|z|^{-\alpha},\lambda_h(q)^{-1})\\
&\leq c\min(|z|^{-\alpha},1),\quad \forall z\in\Gamma_{\theta,\sigma}^\tau.
\end{split}
\end{equation}

Now we can give smoothing properties of the solution operators $F_{h,\tau}^n(q)$ and $E_{h,\tau}^n(q)$.
The proof of the cases $s=0$ and $s=1$ can be found in \cite[Lemma 4.3]{zhang2022identification},
and the case $0<s<1$ follows from standard interpolation theory.
\begin{lemma}\label{lem:smoothing-disc}
For any $q\in \mathcal{A}$ and any $s\in[0,1]$, there exists $c$, independent of $h$, $\tau$, $t_n$ and $q$, such that
for all $v_h\in X_h$
$$t_n^{s\alpha}\|A_h(q)^{s}F_{h,\tau}^n(q)v_h\|_{L^2(\Omega)} + t_n^{1-(1-s)\alpha}\|A_h(q)^{s}E_{h,\tau}^n(q)v_h\|_{L^2(\Omega)}\leq c\|v_h\|_{L^2(\Omega)}.$$
\end{lemma}

The analysis uses frequently the following discrete $L^p(\Omega)$ resolvent estimate.
Note that the estimate \eqref{ineq: Lp dis resol} is different from that in
Lemma \ref{lem:maximum-resolv} in that the latter allows also mappings to
the space $W^{1,\infty}(\Omega)$. This difference has important consequences in the error analysis, and it has to be overcome alternatively.
\begin{lemma}
Let $\Omega$ be a smooth domain or a convex polygon and
$q\in W^{1,\infty}(\Omega)$. Then for any $p\in[1,\infty]$ and $v_h\in X_h$, there holds
\begin{equation}\label{ineq: Lp dis resol}
(1+|z|)\|(z+A_h(q))^{-1}v_h\|_{L^p(\Omega)}\leq c\|v_h\|_{L^p(\Omega)},\quad\forall z\in\Sigma_{\theta} , \ \theta\in(\tfrac\pi2,\pi).
\end{equation}
\end{lemma}
\begin{proof}
See \cite[Theorem 1.1]{crouzeix1994resolvent} for the case $d=1$,
and \cite[Theorem 1.1]{li2017maximal} for the case $d=2,3$.
Note that for $d=2,3$, the work \cite{li2017maximal} discussed only the case $p=\infty$. The case $p = 1$ follows by a duality argument
\begin{align*}
    \|(z+A_h(q))^{-1}v_h\|_{L^1(\Omega)} &= \sup_{\|w\|_{L^\infty(\Omega)}=1}\big((z+A_h(q))^{-1}v_h,w\big) \\
    &= \sup_{\|w\|_{L^\infty(\Omega)}=1}\big(v_h,(z+A_h(q))^{-1}P_hw\big) \\
   & \leq \sup_{\|w\|_{L^\infty(\Omega)}=1}\|v_h\|_{L^1(\Omega)}\|(z+A_h(q))^{-1}P_hw\|_{L^\infty(\Omega)}\\
   &\leq c(1+|z|)^{-1}\|v_h\|_{L^1(\Omega)}.
\end{align*}
The intermediate case $p\in(1,\infty)$ follows from Riesz--Thorin interpolation theorem.
\end{proof}

Now we can give an error bound on the FEM approximation
$U_h^N(\Pi_hq^\dag)$ by problem \eqref{eqn:fem} with $q_h=\Pi_hq^\dag$. This estimate plays a central role in establishing Proposition \ref{prop:err:uT-uhN-and-gqh}.
\begin{lemma}\label{lem:err:uq-uPiq}
If Assumption \ref{ass:disc} holds, then there exists $c>0$, independent of $\tau$, $h$ and $t_n$, such that
\begin{equation*}
	\|u(q^\dag)(t_n)-U_h^n(\Pi_hq^\dag)\|_{L^2(\Omega)}\leq c( t_n^{\alpha-1}\tau+\max(1,t_n^{-\alpha})h^2), \quad n=1,2,\dots,N.
\end{equation*}
\end{lemma}
\begin{proof}
Let $u^n\equiv u(q^\dag)(t_n)$ and $U_h^n\equiv U_h^n(q^\dag)$. Then the following estimate holds
\cite[Lemma A.1]{JinZhou:2023IMA}
\begin{equation*}
	\|u^n-U_h^n\|_{L^2(\Omega)}\leq c(\tau t_n^{\alpha-1}+h^2), \quad n=1,2,\dots,N.
\end{equation*}
Next we prove
\begin{equation*}
	\|U_h^n-U_h^n(\Pi_hq^\dag)\|_{L^2(\Omega)}\leq c\max(1,t_n^{-\alpha})h^2,\quad n=1,2,\dots,N.
\end{equation*}
Using the operator $F_{h,\tau}^n(q)$ in \eqref{equ: dis oper F and E}, $U_h^n$ and
$U_h^n(\Pi_hq^\dag)$ can be represented respectively by
\begin{align*}
	U_h^n&
= F_{h,\tau}^n(q^\dag)U_h^0+ (I-F_{h,\tau}^n(q^\dag))A_h(q^\dag)^{-1}P_hf,\\
	U_h^n(\Pi_hq^\dag)&=F_{h,\tau}^n(\Pi_hq^\dag)U_h^0+ (I-F_{h,\tau}^n(\Pi_hq^\dag))A_h(\Pi_hq^\dag)^{-1}P_hf.
\end{align*}
Thus the error $e_h^n:=U_h^n-U_h^n(\Pi_hq^\dag)$ satisfies $e_h^0=0$ and for $n=1,2,\dots,N$
\begin{align*}
	e_h^n=&(F_{h,\tau}^n(q^\dag)-F_{h,\tau}^n(\Pi_hq^\dag))U_h^0+ (A_h(q^\dag)^{-1}-A_h(\Pi_hq^\dag)^{-1})P_hf \\&+ (-F_{h,\tau}^n(q^\dag)A_h(q^\dag)^{-1}+F_{h,\tau}^n(\Pi_hq^\dag)A_h(\Pi_hq^\dag)^{-1})P_hf:={\rm I_1}+{\rm I_2}+{\rm I_3}.
\end{align*}
It remains to bound the three terms separately. For the term ${\rm I_1}$, by the definition of $F_{h,\tau}^n(q^\dag)$
(with $\sigma=t_n^{-1}$ in $\Gamma^\tau_{\theta,\sigma}$), we have
\begin{align*}
  {\rm I}_1 =  \frac{1}{2\pi\rm i}\int_{\Gamma^\tau_{\theta,\sigma}}\!\!\!\!\!\!\!\!e^{zt_{n-1}}\delta_\tau(e^{-z\tau})^{\alpha-1}\big(\big(\delta_\tau(e^{-z\tau})^{\alpha}+A_h(q^\dag)\big)^{-1}
  -\big(\delta_\tau(e^{-z\tau})^{\alpha}+A_h(\Pi_hq^\dag)\big)^{-1}\big)U_h^0\ \mathrm{d}z.
\end{align*}
It follows directly from the identity $B_1^{-1}-B_2^{-1}=B_1^{-1}(B_2-B_1)B_2^{-1}$ that
\begin{align*}
  & \|\big(\delta_\tau(e^{-z\tau})^{\alpha}+A_h(q^\dag)\big)^{-1}
  -\big(\delta_\tau(e^{-z\tau})^{\alpha}+A_h(\Pi_hq^\dag)\big)^{-1}\|_{L^p(\Omega)\to L^2(\Omega)}\\
=&\|\big(\delta_\tau(e^{-z\tau})^{\alpha}+A_h(q^\dag)\big)^{-1}\big(A_h(\Pi_hq^\dag)-A_h(q^\dag)\big) \big(\delta_\tau(e^{-z\tau})^{\alpha}+A_h(\Pi_hq^\dag)\big)^{-1}\|_{L^p(\Omega)\to L^2(\Omega)} \\
\leq&\big\|\big(\delta_\tau(e^{-z\tau})^{\alpha}+A_h(q^\dag)\big)^{-1}A_h(q^\dag)\big\|_{L^2(\Omega)\to L^2(\Omega)}\big\|A_h(\Pi_hq^\dag)^{-1}-A_h(q^\dag)^{-1}\big\|_{L^p(\Omega)\to L^2(\Omega)}\\
&\qquad \qquad \times \big\|\big(\delta_\tau(e^{-z\tau})^{\alpha}+A_h(\Pi_h q^\dag)\big)^{-1}A_h(\Pi_hq^\dag)\big\|_{L^p(\Omega)\to L^p(\Omega)}.
\end{align*}
Upon recalling the estimates \eqref{inequ:Pi_h approx-2}--\eqref{inequ:Pi_h approx-inf} and using the argument of
\cite[Lemma A.1]{JinZhou:2021SINUM}, we obtain
\begin{align}\label{eqn:bdd-fem}
 \|A_h(\Pi_hq^\dag)^{-1}-A_h(q^\dag)^{-1}\|_{L^p(\Omega)\to L^2(\Omega)}&\leq ch^2,\quad \mbox{with }p>\max(d+\epsilon, 2).
\end{align}
Meanwhile, the discrete resolvent estimate \eqref{ineq: Lp dis resol} implies
\begin{align}\label{eqn:resol-op-bdd}
  \|(\delta_\tau(e^{-z\tau})^{\alpha}+A_h( q) )^{-1}A_h(q)\|_{L^p(\Omega)\to L^p(\Omega)} & \leq c, \quad \forall q\in \mathcal{A}.
\end{align}
Combining the estimates \eqref{eqn:bdd-fem} and \eqref{eqn:resol-op-bdd} gives
\begin{align*}
 \|\big(\delta_\tau(e^{-z\tau})^{\alpha}+A_h(q^\dag)\big)^{-1}
  -\big(\delta_\tau(e^{-z\tau})^{\alpha}+A_h(\Pi_hq^\dag)\big)^{-1}\|_{L^p(\Omega)\to L^2(\Omega)}\leq ch^2.
\end{align*}
Consequently, the following bound on ${\rm I}_1$ holds:
\begin{align*}
\|{\rm I_1}\|_{L^2(\Omega)} &\leq ch^{2}\|u_0\|_{L^p(\Omega)}\int_{\Gamma^\tau_{\theta,\sigma}}|e^{zt_n}||z|^{\alpha-1}\ |\mathrm{d}z|\leq ct_n^{-\alpha}h^2.
\end{align*}
For the term ${\rm I_2}$, by the estimate \eqref{eqn:bdd-fem} and the $L^p(\Omega)$ stability of $P_h$, cf. \eqref{inequ:P_h approx}, we obtain
\begin{equation*}
	\|{\rm I_2}\|_{L^2(\Omega)}\leq\|A_h(q^\dag)^{-1}-A_h(\Pi_hq^\dag)^{-1}\|_{L^p(\Omega)\to L^2(\Omega)}\|P_hf\|_{L^p(\Omega)}\leq ch^2.
\end{equation*}
Last, for the term ${\rm I}_3$, we use the splitting
\begin{align*}
	&-F_{h,\tau}^n(q^\dag)A_h(q^\dag)^{-1}+F_{h,\tau}^n(\Pi_hq^\dag)A_h(\Pi_hq^\dag)^{-1} \\
	=& (F_{h,\tau}^n(\Pi_hq^\dag)-F_{h,\tau}^n(q^\dag))A_h(q^\dag)^{-1} - F_{h,\tau}^n(\Pi_hq^\dag)(A_h(q^\dag)^{-1}-A_h(\Pi_hq^\dag)^{-1}) = {\rm I_{3}^1} + {\rm I_{3}^2}.
\end{align*}
Then the argument for the term ${\rm I}_1$ and the boundedness of the operator $A_h(q^\dag)^{-1}P_h$ in $L^p(\Omega)$ \cite[Section 8.5]{BrennerScott:book2008} imply
$$ \| {\rm I_{3}^1} \|_{L^2(\Omega)} \le c h^2 t_n^{-\alpha} \|A_h(q^\dag)^{-1} P_h f  \|_{L^p(\Omega)} \le c h^2 t_n^{-\alpha} \|f\|_{L^p(\Omega)}. $$
Meanwhile, Lemma \ref{lem:smoothing-disc}, the $L^p(\Omega)$ stability of $P_h$
and the estimate \eqref{eqn:bdd-fem} lead to
\begin{align*}
	\|{\rm I_{3}^2}\|_{L^2(\Omega)}
	&\leq  c\min(1,t_n^{-\alpha}) \|(A_h(\Pi_hq^\dag)^{-1}-A_h(q^\dag)^{-1})P_h f\|_{L^2(\Omega)}\\
	&\leq c\min(1,t_n^{-\alpha}) h^2 \| P_h f\|_{L^p(\Omega)} \leq c\min(1,t_n^{-\alpha}) h^2 \|  f\|_{L^p(\Omega)}.
\end{align*}
The desired estimate now follows by combining the preceding estimates.
\end{proof}

Next we provide a crucial \textit{a priori} estimate of $u(q^\dag)(T)-U_h^N(q^*_h)$ and $\nabla q_h^*$.
\begin{proposition}\label{prop:err:uT-uhN-and-gqh}
Let Assumption \ref{ass:disc} hold, and $q^*_h$ be a minimizer of problem \eqref{eqn:obj-disc}-\eqref{eqn:fem}. Then there exists $c$, independent of $\tau$, $h$, $\delta$, $\gamma$ and $T$, such that
\begin{equation*}
\|u(q^\dag)(T)-U_h^N(q^*_h)\|_{L^2(\Omega)}+\gamma^\frac12\|\nabla q_h^*\|_{L^2(\Omega)}\leq c(T^{\alpha-1}\tau+\max(1,T^{-\alpha})h^2+\delta+\gamma^\frac12).
\end{equation*}
\end{proposition}
\begin{proof}
Let $u\equiv u(q^\dag)$. Since $q_h^*$ minimizes problem \eqref{eqn:obj-disc}-\eqref{eqn:fem} and $\Pi_hq^\dagger\in \mathcal{A}_h$, we have
\begin{equation*}
	J_{\gamma,h,\tau}(q_h^*)\leq  J_{\gamma,h,\tau}(\Pi_hq^\dagger).
\end{equation*}
By the $H^1(\Omega)$ stability of the operator $\Pi_h$, cf. \eqref{inequ:Pi_h approx-2},
and Lemma \ref{lem:err:uq-uPiq}, we obtain
\begin{align*}
   &\, \|U_h^N(q^*_h)-z^\delta\|^2_{L^2(\Omega)}+\gamma\|\nabla q_h^*\|^2_{L^2(\Omega)}\\
  \leq  & \, \|U_h^N(\Pi_hq^\dagger)-z^\delta\|^2_{L^2(\Omega)}+\gamma\|\nabla \Pi_hq^\dagger\|^2_{L^2(\Omega)} \\
  \leq  & \, c\big(\|U_h^N(\Pi_hq^\dagger)-u(T)\|^2_{L^2(\Omega)}+\|u(T)-z^\delta\|^2_{L^2(\Omega)}+\gamma\big)\\
  \leq & \, c(T^{2\alpha-2}\tau^2+\max(1,T^{-2\alpha})h^4+\delta^2+\gamma).
\end{align*}
Then the triangle inequality yields
\begin{align*}
	&\|u(T)-U^N_h(q_h^*)\|^2_{L^2(\Omega)}+\gamma\|\nabla q_h^*\|^2_{L^2(\Omega)} \\
\leq \, & c\big(\|u(T)-z^\delta\|^2_{L^2(\Omega)}+\|z^\delta-U^N_h(q_h^*)\|^2_{L^2(\Omega)}+\gamma\|\nabla q_h^*\|^2_{L^2(\Omega)}\big) \\
\leq \, & c(T^{2\alpha-2}\tau^2+\max(1,T^{-2\alpha})h^4+\delta^2+\gamma).
\end{align*}
This completes the proof of the lemma.
\end{proof}

\subsection{Bound on the error $\bar{\partial}_\tau^\alpha U_h^n(q^*_h)-\partial_t^\alpha u(q^\dag)(t_n)$}
Next, we estimate the decay of the discrete (fractional) derivative $ \bar{\partial}_\tau^\alpha U_h^n(q^*_h)$ and bound the term
$\bar{\partial}_\tau^\alpha U_h^n(q^*_h)-\partial_t^\alpha u(q^\dag)(t_n)$ in terms of
$\|q_h^*-q^\dag\|_{L^2(\Omega)}$; see Lemma \ref{lem:est-dt-uhn} and Proposition \ref{prop:err-dtauuq-uqh} for the precise statement.
These estimates play a central role in establishing Theorem \ref{thm:main-disc}.
We need the following time semidiscrete scheme for problem \eqref{eqn:fde}: Find $U^n\equiv U^n(q)\in H_0^1(\Omega)$
with $U^0 = u_0$ such that
\begin{equation}\label{equ: time stepping scheme}
	\bar{\partial}_\tau^\alpha U^n(q) + A(q)U^n(q) = f, \quad n=1,2,\dots,N.
\end{equation}
The discrete Laplace transform gives
\begin{equation}\label{eqn:sol-rep-semi}
	U^n=F^n_\tau(q)u_0 + \tau\sum_{j=1}^nE^j_\tau(q)f=F^n_\tau(q)u_0 + (I-F^n_\tau(q))A(q)^{-1}f,
\end{equation}
where the time-semidiscrete solution operators $F_\tau^n\tau(q)$ and $E_\tau^n(q)$ are defined respectively by \cite[Section 3.2]{JinZhou:2023book}
\begin{equation}\label{equ: time stepping dis oper F and E}
\begin{aligned}
F^n_\tau(q)&=\frac{1}{2\pi i}\int_{\Gamma^\tau_{\theta,\sigma}}e^{zt_{n-1}}\delta_\tau(e^{-z\tau})^{\alpha-1}(\delta_\tau(e^{-z\tau})^{\alpha}+A(q))^{-1}\ \mathrm{d}z,\\
E^n_\tau(q)&=\frac{1}{2\pi i}\int_{\Gamma^\tau_{\theta,\sigma}}e^{zt_{n-1}}(\delta_\tau(e^{-z\tau})^{\alpha}+A(q))^{-1}\ \mathrm{d}z.
\end{aligned}
\end{equation}

The next lemma gives a temporal error estimate for the approximate time (fractional) derivative.
\begin{lemma}\label{lem: err on frac discretization}
Let Assumption \ref{ass:disc} hold, and $u(q^\dag)$ and $\{U^n(q^\dag)\}_{n=0}^N$ be the solutions of problems \eqref{eqn:fde} and \eqref{equ: time stepping scheme} for $q^\dag$, respectively. Then there exists $c>0$, independent of $\tau$, $t_n$ and $q^\dag$, such that
	\begin{equation*}
		\|\partial_t^\alpha u(q^\dag)(t_n)-\bar{\partial}_\tau^\alpha U^n(q^\dag)\|_{L^2(\Omega)}\leq c\tau t_n^{-\alpha-1}.
	\end{equation*}
\end{lemma}
\begin{proof}
Let $u\equiv u(q^\dag)$ and $U^n\equiv U^n(q^\dag)$. Then $W^n :=\bar{\partial}_\tau^\alpha U^n$ satisfies $W^0=f-A(q^\dag)u_0$ and
\begin{equation*}
	\bar{\partial}_\tau^\alpha W^n +A(q^\dag)W^n=0,\quad n=1,2,\ldots,N.
\end{equation*}
It follows from the solution representations \eqref{eqn:sol} and \eqref{eqn:sol-rep-semi} that
\begin{equation*}
	\partial_t^\alpha u(t_n) = F(t_n;q^\dag)(f-A(q^\dag)u_0) \quad \mbox{and} \quad
    \bar{\partial}_\tau^\alpha U^n = F^n_\tau(q^\dag)(f-A(q^\dag)u_0).
\end{equation*}
It follows from the estimate 
$\|F(t_n;q)-F^n_\tau(q)\|_{L^2(\Omega)\to L^2(\Omega)}\leq  cn^{-1}t_n^{-\alpha}$ \cite[Lemma 15.6]{JinZhou:2023book}
and Assumption \ref{ass:disc} that
\begin{align*}
	\|\partial_t^\alpha u(q^\dag)(t_n)-\bar{\partial}_\tau^\alpha U^n(q^\dag)\|_{L^2(\Omega)} &\leq \|F(t_n;q^\dag)-F^n_\tau(q^\dag)\|_{L^2(\Omega)\to L^2(\Omega)}\|f-A(q^\dag)u_0\|_{L^2(\Omega)}\\&\leq  c \tau t_n^{-\alpha-1}(\|f\|_{L^2(\Omega)}+\|u_0\|_{H^2(\Omega)})\leq c \tau t_n^{-\alpha-1}.
\end{align*}
This completes the proof of the lemma.
\end{proof}

The next lemma bounds the error between the discrete (fractional) derivative due to spatial discretization.
\begin{lemma}\label{lem: err on space discretization}
Let Assumption \ref{ass:disc} hold, and $\{U^n(q^\dag)\}_{n=0}^N$ and $\{U^n_h(q^\dag)\}_{n=0}^N$ be the solutions of problem \eqref{equ: time stepping scheme} with $q^\dag$ and problem \eqref{eqn:fem} with $q^\dag$, respectively. Then there exists $c>0$, independent of $\tau$, $h$, $t_n$ and $q^\dag$, such that
\begin{equation*}
	\|\bar{\partial}_\tau^\alpha \big(U^n(q^\dag)- U_h^n(q^\dag)\big)\|_{L^2(\Omega)}\leq  ch^2\max(t_n^{-\alpha},t_n^{-2\alpha}).
\end{equation*}
Moreover, if $u_0\in W^{2,\infty}(\Omega)$, then with $\ell_h:=|\log h|$,
\begin{equation*}
   \|\bar{\partial}_\tau^\alpha \big(U^n(q^\dag)-U_h^n(q^\dag)\big)\|_{L^\infty(\Omega)}\leq  cht_n^{-\frac\alpha2} + ch^2\ell_h^2\max(t_n^{-\alpha}, t_n^{-2\alpha}).
\end{equation*}
\end{lemma}
\begin{proof}
Let $U^n\equiv U^n(q^\dag)$ and $U^n_h\equiv U^n_h(q^\dag)$, and $e^n_{h,\tau}:=\bar{\partial}_\tau^\alpha(U^n - U_h^n)$.
To bound $e^n_{h,\tau}$ in $L^2(\Omega)$, the solution representation \eqref{eqn:sol-rep-semi} yields
\begin{equation}\label{equ: dt time stepping error}
\begin{aligned}
e^n_{h,\tau} & = F^n_\tau(q^\dag)(f-A(q^\dag)u_0)-F^n_{h,\tau}(q^\dag)(P_hf-A_h(q^\dag)P_hu_0) \\
& = \big(F^n_\tau(q^\dag)-F^n_{h,\tau}(q^\dag)P_h\big)f + \big(F_{h,\tau}^n(q^\dag)A_h(q^\dag)P_h-F^n_\tau(q^\dag)A(q^\dag)\big)u_0\\
  &:= {\rm I}_1+{\rm I}_2.
\end{aligned}
\end{equation}
Now we bound the two terms ${\rm I}_1$ and ${\rm I}_2$ separately. Let $B_{h,\tau}=\big(\delta_\tau(e^{-z\tau})^{\alpha}+A(q^\dag)\big)^{-1}-\big(\delta_\tau(e^{-z\tau})^{\alpha}+A_h(q^\dag)\big)^{-1}P_h$. It follows from the estimates in \eqref{ineq: dis kernel property} that for all $z\in\Gamma_{\theta,\sigma}^\tau$ \cite[Theorem 5.2 and Remark 7.4]{fujita1991evolution}
\begin{equation}\label{eqn:esti-res-app}
  \|B_{h,\tau}\|_{L^2(\Omega)\to L^2(\Omega)}\leq ch^2.
\end{equation}
Then choosing $\sigma=t_n^{-1}$ in the contour $\Gamma_{\theta,\sigma}^\tau$ leads to
    \begin{align*}
    	\|{\rm I}_1\|_{L^2(\Omega)}
    	&\leq c\|f\|_{L^2(\Omega)} \int_{\Gamma^\tau_{\theta,\sigma}}|e^{zt_n}||\delta_\tau(e^{-z\tau})|^{\alpha-1}\|B_{h,\tau}\|_{L^2(\Omega)\to L^2(\Omega)}\ |\mathrm{d}z| \\
    	&\leq ch^2\|f\|_{L^2(\Omega)}\int_{\Gamma^\tau_{\theta,\sigma}}|e^{zt_n}||z|^{\alpha-1}\ |\mathrm{d}z|\leq ch^2t_n^{-\alpha}.
    \end{align*}
To estimate the term ${\rm I}_2$, by the identity
$\big(\delta_\tau(e^{-z\tau})^{\alpha}+A_h(q^\dag)\big)^{-1}A_h(q^\dag)P_h-\big(\delta_\tau(e^{-z\tau})^{\alpha}+A(q^\dag)\big)^{-1}A(q^\dag)= \big(P_h-I\big) +  \delta_\tau(e^{-z\tau})^{\alpha}B_{h,\tau},$
we derive
\begin{align*}
{\rm I}_2 &= \frac{1}{2\pi\mathrm{i}}\int_{\Gamma_{\theta,\sigma}^\tau} e^{zt_n} \delta_\tau(e^{-z\tau})^{\alpha-1}(P_h - I) u_0 \,\mathrm{d} z + \frac{1}{2\pi\mathrm{i}}\int_{\Gamma_{\theta,\sigma}^\tau} e^{zt_n} \delta_\tau(e^{-z\tau})^{2\alpha-1}B_{h,\tau}u_0{\rm d}z.
\end{align*}
Then with $\sigma=t_n^{-1}$ in the contour $\Gamma_{\theta,\sigma}^\tau$, by the estimates \eqref{eqn:esti-res-app},
\eqref{ineq: dis kernel property}, and \eqref{inequ:P_h approx}, we derive
\begin{align*}
	\|{\rm I}_2\|_{L^2(\Omega)}	&\leq ch^2 \int_{\Gamma^\tau_{\theta,\sigma}}|e^{zt_n}|(|z|^{\alpha-1}\|u_0\|_{H^2(\Omega)} + |z|^{2\alpha-1}\|u_0\|_{L^2(\Omega)})\ |\mathrm{d}z| \\
    &\leq ch^2t_n^{-\alpha}\|u_0\|_{H^2(\Omega)}+ch^2t_n^{-2\alpha}\|u_0\|_{L^2(\Omega)}\leq ch^2\max(t_n^{-\alpha},t_n^{-2\alpha}).
\end{align*}
To bound $\|e^n_{h,\tau}\|_{L^\infty(\Omega)}$, we split $e^n_{h,\tau}$ into
    \begin{equation*}
    	e^n_{h,\tau}:=\big(\bar{\partial}_\tau^\alpha U^n - P_h\bar{\partial}_\tau^\alpha U^n\big)+ \big(P_h\bar{\partial}_\tau^\alpha U^n - \bar{\partial}_\tau^\alpha U_h^n\big):= {\rm I}_3+{\rm I}_4.
    \end{equation*}
It follows from the estimates in \eqref{ineq: dis kernel property}, \eqref{ineq:max-resolvents}, and the argument of Lemma \ref{lem:dut} that
$$\|F^n_\tau(q^\dag)\|_{L^\infty(\Omega)\to W^{1,\infty}(\Omega)}\leq ct_n^{-\frac\alpha2}.$$
From the approximation property of $P_h$ in \eqref{inequ:P_h approx} and the assumption $u_0\in W^{2,\infty}(\Omega)$, we have
\begin{align*}
  &\|{\rm I}_3\|_{L^\infty(\Omega)}  \leq ch\|\bar{\partial}_\tau^\alpha U^n\|_{W^{1,\infty}(\Omega)} \\
  \leq& ch\|F^n_\tau(q^\dag)\|_{L^\infty(\Omega)\to W^{1,\infty}(\Omega)}\big(\|f\|_{L^\infty(\Omega)} +\|A(q^\dag)u_0\|_{L^\infty(\Omega)}\big)\leq cht_n^{-\frac\alpha2}.
\end{align*}
It remains to bound the term ${\rm I}_4$. From the representation \eqref{equ: dt time stepping error}, we obtain
\begin{equation*}
   {\rm I}_4 = \big(P_hF^n_\tau(q^\dag)-F_{h,\tau}^n(q^\dag)P_h\big)f + \big(F_{h,\tau}^n(q^\dag)A_h(q^\dag)P_h-P_hF^n_\tau(q^\dag)A(q^\dag)\big)u_0.
\end{equation*}
Let $R_h$ be the standard Ritz projection. Then direct computation gives
\begin{align*}
   &P_h\big(\delta_\tau(e^{-z\tau})^{\alpha}+A(q^\dag)\big)^{-1} - \big(\delta_\tau(e^{-z\tau})^{\alpha}+A_h(q^\dag)\big)^{-1}P_h\\
    = &A_h(q^\dag)\big(\delta_\tau(e^{-z\tau})^{\alpha}+A_h(q^\dag)\big)^{-1}(P_h-R_h)\big(\delta_\tau(e^{-z\tau})^{\alpha}+A(q^\dag)\big)^{-1}.
\end{align*}
Hence, we obtain
\begin{align*}
&P_hF^n_\tau(q^\dag)-F_{h,\tau}^n(q^\dag)P_h =   \frac{1}{2\pi\rm i}\int_{\Gamma^\tau_{\theta,\sigma}}e^{zt_{n-1}}\delta_\tau(e^{-z\tau})^{\alpha-1} \\
&\qquad\qquad \times \big[ A_h(q^\dag)\big(\delta_\tau(e^{-z\tau})^{\alpha}+A_h(q^\dag)\big)^{-1}(P_h-R_h)\big(\delta_\tau(e^{-z\tau})^{\alpha}+A(q^\dag)\big)^{-1}\big]\ \mathrm{d}z
\end{align*}
and
\begin{align*}
&F_{h,\tau}^n(q^\dag)A_h(q^\dag)P_h-P_hF^n_\tau(q^\dag)A(q^\dag) =\frac{1}{2\pi\rm i}\int_{\Gamma^\tau_{\theta,\sigma}}e^{zt_{n-1}}\delta_\tau(e^{-z\tau})^{2\alpha-1} \\
&\qquad \qquad \times \big[A_h(q^\dag)\big(\delta_\tau(e^{-z\tau})^{\alpha}+A_h(q^\dag)\big)^{-1}(P_h-R_h)\big(\delta_\tau(e^{-z\tau})^{\alpha}+A(q^\dag)\big)^{-1}\big]\ \mathrm{d}z.
\end{align*}
The resolvent estimates \eqref{ineq: Lp dis resol} and \eqref{ineq:max-resolvents} imply
\begin{align*}
 & \|A_h(q^\dag)\big(\delta_\tau(e^{-z\tau})^{\alpha}+A_h(q^\dag)\big)^{-1}\|_{L^\infty(\Omega)\to L^\infty(\Omega)}\leq c,\\
 & \|A(q^\dag)\big(\delta_\tau(e^{-z\tau})^{\alpha}+A(q^\dag)\big)^{-1}\|_{L^\infty(\Omega)\to L^\infty(\Omega)}\leq c.
\end{align*}
Now the estimate 
$\|(P_h-R_h)A(q^\dag)^{-1}\|_{L^\infty(\Omega)\to L^\infty(\Omega)}\leq ch^2\ell_h^2$ holds (\cite[p. 1658]{palencia1996maximum} and \cite[p. 220]{bakaev2001maximum}).
Then letting $\sigma=t_n^{-1}$ in the contour $\Gamma_{\theta,\sigma}^\tau$ leads to
\begin{equation*}
\|P_hF^n_\tau(q^\dag)-F_{h,\tau}^n(q^\dag)P_h\|_{L^\infty(\Omega)\to L^\infty(\Omega)}\leq ch^2\ell_h^2t_n^{-\alpha},
\end{equation*}
and
\begin{equation*}
\|F_{h,\tau}^n(q^\dag)A_h(q^\dag)P_h-P_hF^n_\tau(q^\dag)A(q^\dag)\|_{L^\infty(\Omega)  \to L^\infty(\Omega)}\leq ch^2\ell_h^2 t_n^{-2\alpha},
\end{equation*}
which implies
\begin{equation*}
\|{\rm I}_4\|_{L^\infty(\Omega)} \leq ch^2\ell_h^2t_n^{-\alpha}\|f\|_{L^\infty(\Omega)}+ch^2\ell_h^2t_n^{-2\alpha}\|u_0\|_{L^\infty(\Omega)}\leq ch^2\ell_h^2\max(t_n^{-\alpha}, t_n^{-2\alpha}).
\end{equation*}
Combining the preceding estimates completes the proof.
\end{proof}

The next result gives a discrete analogue of Lemma \ref{lem:dut}.
\begin{lemma}\label{lem:est-dt-uhn}
Let $\Omega$ be a convex $C^2$ domain, and Assumption \ref{ass:disc} hold. Let  $\{U^n_h(q^\dag)\}_{n=0}^N$ be the solution of problem \eqref{eqn:fem} with $q^\dag$. Then for small $h>0$, there exists $c>0$, independent of $\tau$, $h$, $t_n$ and $q^\dag$, such that
	\begin{equation*}
		\|\bar{\partial}_\tau^\alpha U_h^n(q^\dag)\|_{W^{1,\infty}(\Omega)}\leq c\max(t_n^{-\frac\alpha2},t_n^{-\alpha},t_n^{-2\alpha}),\quad  n=1,2,\dots,N.
	\end{equation*}
\end{lemma}
\begin{proof}
It follows from the inverse estimate on the FEM space $X_h$ \cite[equation (1.12)]{Thome2006GalerkinFE}, Lemma \ref{lem: err on space discretization} and
the approximation property and $W^{1,\infty}(\Omega)$ stability of $P_h$ in \eqref{inequ:P_h approx} that
\begin{align*}
		\|\bar{\partial}_\tau^\alpha U_h^n(q^\dag)\|_{W^{1,\infty}(\Omega)}&\leq \|\bar{\partial}_\tau^\alpha U_h^n(q^\dag)-P_h\bar{\partial}_\tau^\alpha U^n(q^\dag)\|_{W^{1,\infty}(\Omega)} + \|P_h\bar{\partial}_\tau^\alpha U^n(q^\dag)\|_{W^{1,\infty}(\Omega)} \\
		&\leq ch^{-1}\|\bar{\partial}_\tau^\alpha U_h^n(q^\dag)-\bar{\partial}_\tau^\alpha U^n(q^\dag)\|_{L^\infty(\Omega)} + c\|\bar{\partial}_\tau^\alpha U^n(q^\dag)\|_{W^{1,\infty}(\Omega)} \\
		&\leq ct_n^{-\frac\alpha2}+ch\ell_h^2\max(t_n^{-\alpha},t_n^{-2\alpha})\leq c \max(t_n^{-\frac\alpha2},t_n^{-\alpha},t_n^{-2\alpha}).
\end{align*}
This completes the proof of the lemma.
\end{proof}

We also have a discrete version of Lemma \ref{lem:dt:u1-u2}. This estimate is
crucial to the error analysis.
\begin{proposition}\label{prop:err-dtauuq-uqh}
Let $\Omega$ be a convex  $C^2$  domain, and Assumption \ref{ass:disc} hold. Let  $\{U^n_h(q^\dag)\}_{n=0}^N$ and $\{U^n_h(q_h^*)\}_{n=0}^N$ be the solutions of problem \eqref{eqn:fem} with $q^\dag$ and $q_h^*$, respectively. Then for small $\epsilon,h>0$, there exists $c>0$, independent of $\tau$, $h$, and $t_n$, such that for all $n=1,\ldots,N$
	\begin{equation*}
		\|\bar{\partial}_\tau^\alpha \big(U_h^n(q^\dag)- U_h^n(q_h^*)\big)\|_{L^2(\Omega)}\leq c\max(t_n^{-\alpha}, t_n^{-\frac\alpha2(1-2\epsilon)},t_n^{-\alpha(1-\epsilon)},t_n^{-\alpha(2-\epsilon)})\|q^\dag-q^*_h\|_{L^2(\Omega)},
	\end{equation*}
\end{proposition}
\begin{proof}
Let $W_h^n:=\bar{\partial}_\tau^\alpha \big(U_h^n(q_h^*)- U_h^n(q^\dag)\big)$. Then $W_h^0=\big(A_h(q^\dag)-A_h(q_h^*)\big)U_h^0$ and
\begin{equation*}
\bar{\partial}_\tau^\alpha W_h^n +A_h(q_h^*)W_h^n= \big(A_h(q^\dag)-A_h(q_h^*)\big)\bar{\partial}_\tau^\alpha  U_h^n(q^\dag),\quad n=1,2,\dots,N.
\end{equation*}
By the discrete Laplace transform, we obtain
\begin{equation*}
    W_h^n= F_{h,\tau}^n(q_h^*)\big(A_h(q^\dag)-A_h(q_h^*)\big)U_h^0 + \tau\sum_{j=1}^n E_{h,\tau}^{n-j}(q_h^*)\big(A_h(q^\dag)-A_h(q_h^*)\big)\bar{\partial}_\tau^\alpha U_h^j(q^\dag).
\end{equation*}
Next we bound the two terms separately. Note that the following inequality holds
$$c_1\|A_h(q_h^*)^\frac12v_h\|_{L^2(\Omega)}\leq \|\nabla v_h\|_{L^2(\Omega)}\leq c_2\|A_h(q_h^*)^{\frac12}v_h\|_{L^2(\Omega)},\quad \forall v_h\in X_h.$$
This, Lemma \ref{lem:smoothing-disc}, the $W^{1,\infty}(\Omega)$ stability of $P_h$ in \eqref{inequ:P_h approx}, and the boundedness of $A_h(q_h^*)^{-\frac12}$ in $L^2(\Omega)$ imply
\begin{equation*}
\begin{split}
&\qquad \, \|F_{h,\tau}^n(q_h^*)\big(A_h(q^\dag)-A_h(q_h^*)\big)U_h^0\|_{L^2(\Omega)}\\ &=\sup_{\|\varphi_h\|_{L^{2}(\Omega)}=1}\big((q^\dag-q_h^*)\nabla P_hu_0,\nabla A_h(q_h^*)^{\frac12}F_{h,\tau}^n(q_h^*)A_h(q_h^*)^{-\frac12}\varphi_h\big) \\
&\leq c\|q^\dag-q_h^*\|_{L^{2}(\Omega)}\|\nabla P_hu_0\|_{L^\infty(\Omega)}\|A_h(q_h^*)F_{h,\tau}^n(q_h^*)\|_{L^{2}(\Omega) \to L^{2}(\Omega)} \\
& \qquad \qquad \qquad\qquad \qquad \qquad  \times \sup_{\|\varphi_h\|_{L^{2}(\Omega)}=1}\|A_h(q_h^*)^{-\frac12}\varphi_h\|_{L^{2}(\Omega)}\\
&\leq ct_n^{-\alpha}\|q^\dag-q_h^*\|_{L^{2}(\Omega)}.
\end{split}
\end{equation*}
Similarly, the boundedness of the operator $A(q_h^*)^{-\frac12+\epsilon}$ in $L^2(\Omega)$ and Lemmas \ref{lem:smoothing-disc}-\ref{lem:est-dt-uhn} lead to
\begin{equation*}
\begin{split}
&\qquad\,\|E_{h,\tau}^{n-j}(q_h^*)\big(A_h(q^\dag)-A_h(q_h^*)\big)\bar{\partial}_\tau^\alpha U_h^j(q^\dag)\|_{L^2(\Omega)}\\ &=\sup_{\|\varphi_h\|_{L^{2}(\Omega)}=1}\big((q^\dag-q_h^*)\nabla \bar{\partial}_\tau^\alpha  U_h^j(q^\dag),\nabla A_h(q_h^*)^{\frac12-\epsilon}E_{h,\tau}^{n-j}(q_h^*)A_h(q_h^*)^{-\frac12+\epsilon}\varphi_h\big) \\
&\leq c\|q^\dag-q_h^*\|_{L^{2}(\Omega)}\|\nabla \bar{\partial}_\tau^\alpha  U_h^j(q^\dag)\|_{L^\infty(\Omega)}\|A_h(q_h^*)^{1-\epsilon}E_{h,\tau}^{n-j}(q_h^*)\|_{L^{2}(\Omega) \to L^{2}(\Omega)}\\
&\qquad\qquad\qquad\qquad\qquad\qquad\times \sup_{\|\varphi_h\|_{L^{2}(\Omega)}=1}\|A_h(q_h^*)^{-\frac12+\epsilon}\varphi_h\|_{L^{2}(\Omega)} \\
&\leq ct_{n-j}^{-1+\epsilon\alpha}\max(t_j^{-\frac\alpha2},t_j^{-\alpha},t_j^{-2\alpha})\|q^\dag-q_h^*\|_{L^{2}(\Omega)}.
\end{split}
\end{equation*}
Finally, the preceding two estimates and the triangle inequality imply
    \begin{align*}
    	\|W_h^n\|_{L^2(\Omega)}&\leq ct_n^{-\alpha}\|q^\dag-q_h^*\|_{L^{2}(\Omega)} + c\|q^\dag-q_h^*\|_{L^{2}(\Omega)}\tau\sum_{j=1}^nt_{n-j}^{-1+\epsilon\alpha}\max(t_j^{-\frac\alpha2},t_j^{-\alpha},t_j^{-2\alpha})\\&\leq c\max(t_n^{-\alpha}, t_n^{-\frac\alpha2(1-2\epsilon)},t_n^{-\alpha(1-\epsilon)},t_n^{-\alpha(2-\epsilon)})\|q^\dag-q^*_h\|_{L^2(\Omega)}.
    \end{align*}
This completes the proof of the lemma.
\end{proof}

\section{Numerical experiments and discussions}
\label{sec:numer}
Now we present numerical results for the time fractional diffusion model. We employ the conjugate gradient method \cite{Alifanov:1995} to solve the discrete optimization problem. The gradient $J'_\gamma$ is computed using an alternative adjoint technique  \cite[Section 5]{cheng2020inverse}.
The details of the algorithm are described in Algorithm \ref{alg:cgm} in the appendix for completeness. The lower and upper bounds of the admissible set $\mathcal{A}$ are taken to be $0.5$ and $5.0$. The noise data $z^\delta$ is generated by
\begin{equation*}
	z^\delta(x)=u(q^\dag)(x,T)+\epsilon\|u(q^\dagger)(T)\|_{L^\infty(\Omega)}\xi(x),\quad x\in \Omega,
\end{equation*}
where $\xi$ follows the standard Gaussian distribution and $\epsilon>0$ is the relative noise level. To measure the convergence of the approximation $q_h^*$, we employ two metrics, i.e., $e_q=\|q^\dagger-q^*_h\|_{L^2(\Omega)}$ and $e_u=\|u(q^\dag)(T)-u^N_h(q^*_h)\|_{L^2(\Omega)}$.

First, we consider the one-dimension case with varying terminal time $T$ and fractional order $\alpha$.
\begin{example}\label{exam:1d}
Let $\Omega=(0,1)$, $q^\dag(x)=\max(\min(1+\frac14\sin(\pi x),\frac{319}{256}), \frac{67}{64})$, $u_0(x)=x(1-x)$ and $f\equiv1$. Consider the following two cases: (a) $T=1.00$ and $\alpha=0.25$, $0.50$ and $0.75$; (b) $\alpha=0.50$ and $T=10^{-5}$, $3.00$ and $5.00$.
\end{example}

The exact data is obtained using a finer grid with a mesh size $h=1/1600$ and the number  $N=1280$ of time steps. The numerical results for Example \ref{exam:1d} are presented in Table \ref{tab:err-1d}. Note that when $T$ is large, the theory predicts a convergence rate $O(\delta^\frac14)$ at best for $e_q$, and $O(\delta)$ for $e_u$ of the state approximation (if the parameters are chosen properly). For case (a), both $e_q$ and $e_u$ exhibit a clear decay property as the noise level $\delta$ tends to zero but the empirical rate of $e_q$ is faster than the theoretical one, indicating room for further improvement in convergence analysis. In Fig. \ref{fig:recon-frac1d}, we present the numerical reconstructions of case (a) at different noise levels. The results for case (b) in Table \ref{tab:err-1d} show that the convergence behaviors of both $e_q$ and $e_u$ fail to hold due to the loss of conditional stability in Theorems \ref{thm:main-cont} and \ref{thm:main-disc} when the terminal time $T$ is sufficiently small.

\begin{table}[hbt!]

\centering

\caption{Numerical results for Example \ref{exam:1d}, initialized with $M=113$ and the total time step $N=30$. \label{tab:err-1d}}
\begin{threeparttable}
\subfigure[Results for case (a), with $T$ fixed at $T=1.00$, and varying $\alpha$.]{
\begin{tabular}{c|c|ccccc}
\toprule
$\alpha$
& $\epsilon$  & 1.00e-2 & 5.00e-3 & 2.50e-3 & 1.00e-3 &  rate\\
& $\gamma$    & 4.00e-8 & 1.00e-8 & 2.50e-9 & 4.00e-10&   \\
\midrule
$0.25$ & $e_q$ & 2.67e-2 & 1.76e-2 & 1.54e-2 & 8.42e-3 & 0.48 \\
       & $e_u$ & 2.53e-4 & 1.26e-4 & 7.90e-5 & 3.50e-5 & 0.86\\
\midrule
$0.50$ & $e_q$ & 2.56e-2 & 1.85e-2 & 1.40e-2 & 6.28e-3 & 0.58 \\
       & $e_u$ & 2.77e-4 & 1.17e-4 & 6.51e-5 & 3.29e-5 & 0.94 \\
\midrule
$0.75$ & $e_q$ & 2.57e-2 & 1.72e-2 & 1.46e-2 & 6.37e-3 & 0.57 \\
       & $e_u$ & 2.60e-4 & 1.03e-4 & 7.22e-5 & 3.97e-5 & 0.83 \\
\bottomrule
\end{tabular}}
\subfigure[Results for case (b) with $\alpha$ fixed at $0.5$, and varying $T$.]{
\begin{tabular}{c|c|ccccc}
\toprule
$T$& $\epsilon$  & 1.00e-2 & 5.00e-3 & 2.50e-3 & 1.00e-3 &  trend\\
& $\gamma$    & 4.00e-6 & 1.00e-6 & 2.50e-7 & 4.00e-8&   \\
\midrule
$10^{-5}$   & $e_q$ & 9.06e-2 & 9.93e-2 & 9.82e-2 & 1.00e-1 & ---\\
              & $e_u$ & 8.27e-4 & 7.79e-4 & 7.63e-4 & 7.61e-4 & ---\\
\midrule
$3.00$      & $e_q$ & 4.31e-2 & 2.42e-2 & 1.87e-2 & 1.48e-2 & $\searrow$\\
              & $e_u$ & 5.83e-4 & 2.43e-4 & 1.43e-4 & 6.71e-5 & $\searrow$ \\
              \midrule
$5.00$      & $e_q$ & 3.66e-2 & 2.94e-2 & 1.84e-2 & 1.51e-2 & $\searrow$ \\
              & $e_u$ & 4.46e-4 & 2.97e-4 & 1.16e-4 & 7.05e-5 & $\searrow$\\
\bottomrule
\end{tabular}}
\end{threeparttable}
\end{table}

\begin{figure}[htb!]
\centering
\setlength{\tabcolsep}{0pt}
\begin{tabular}{ccc}
\includegraphics[width=0.33\textwidth,trim={4.5cm 9cm 4cm 9cm}, clip]{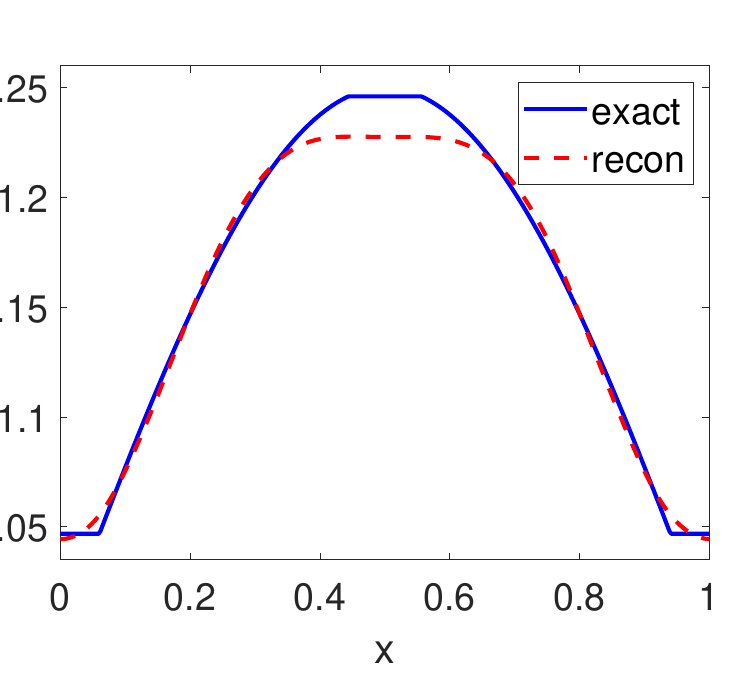} & \includegraphics[width=0.33\textwidth,trim={4.5cm 9cm 4cm 9cm}, clip]{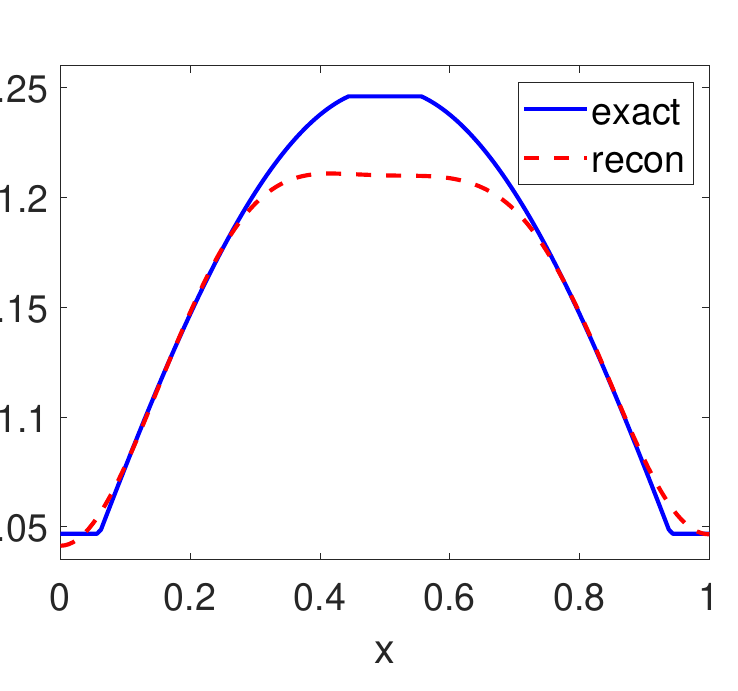} & \includegraphics[width=0.33\textwidth,trim={4.5cm 9cm 4cm 9cm}, clip]{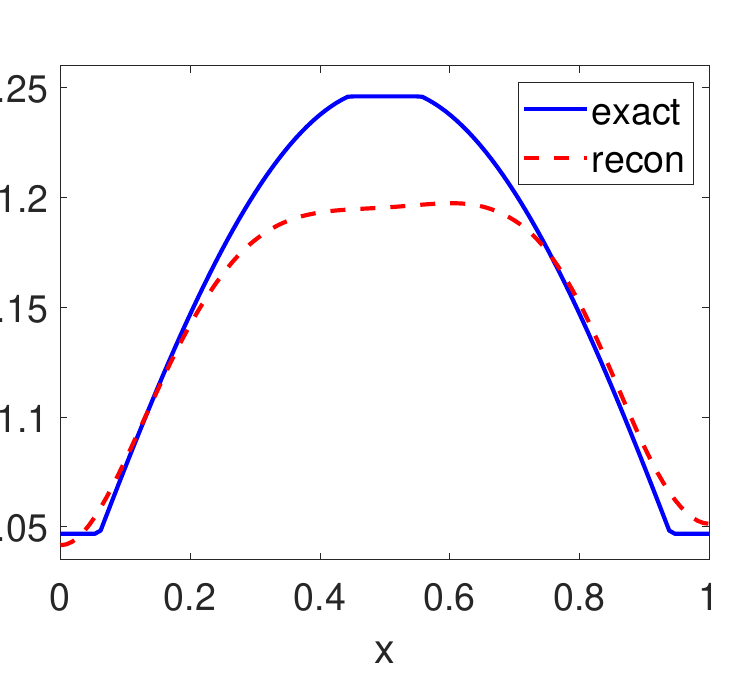}\\
 \includegraphics[width=0.33\textwidth,trim={4.5cm 9cm 4cm 9cm}, clip]{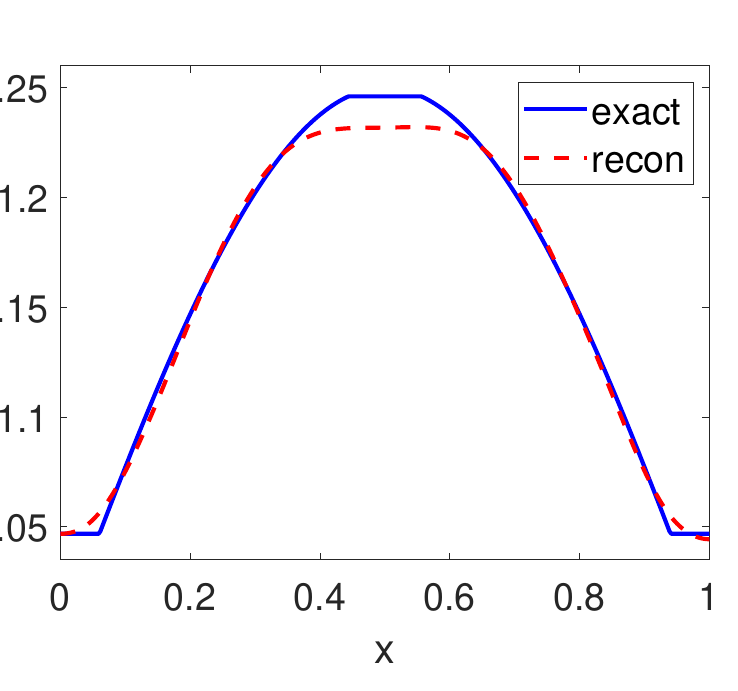} & \includegraphics[width=0.33\textwidth,trim={4.5cm 9cm 4cm 9cm}, clip]{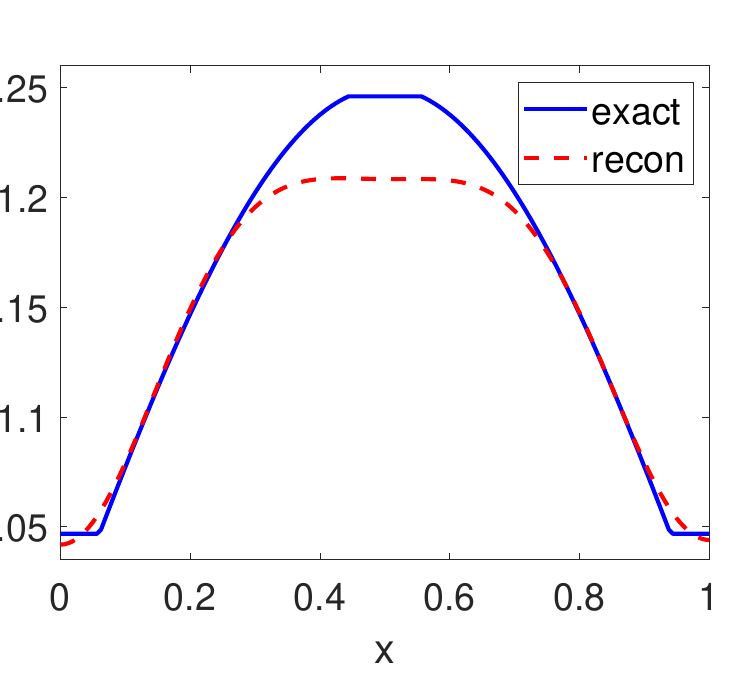} & \includegraphics[width=0.33\textwidth,trim={4.5cm 9cm 4cm 9cm}, clip]{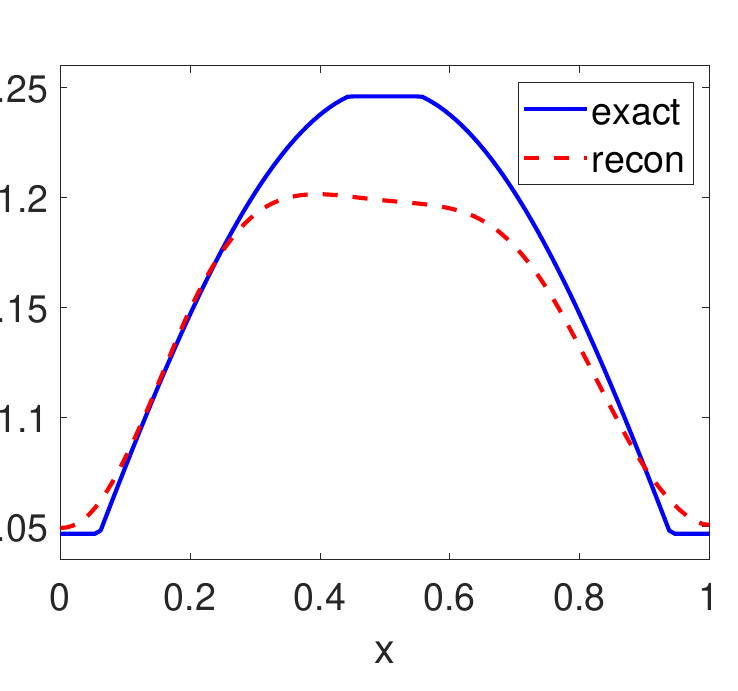}\\
 \includegraphics[width=0.33\textwidth,trim={4.5cm 9cm 4cm 9cm}, clip]{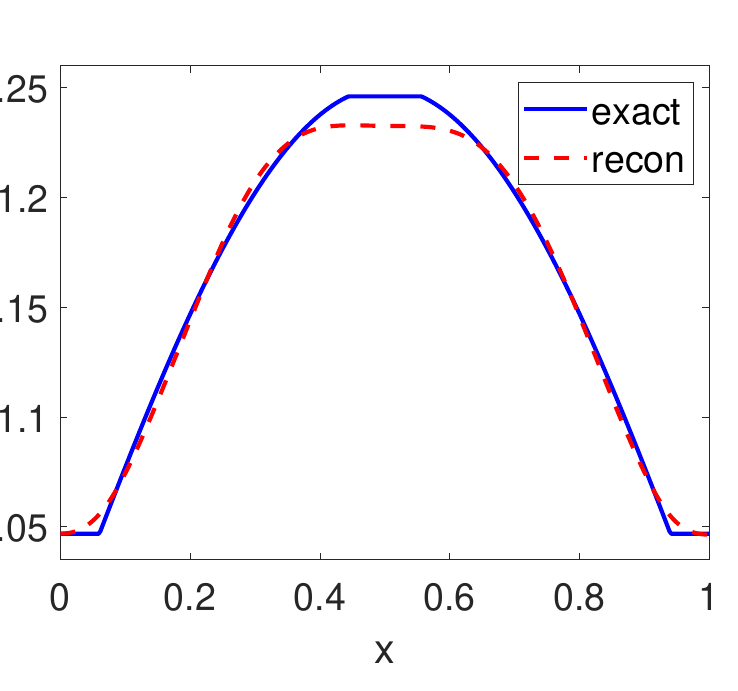} & \includegraphics[width=0.33\textwidth,trim={4.5cm 9cm 4cm 9cm}, clip]{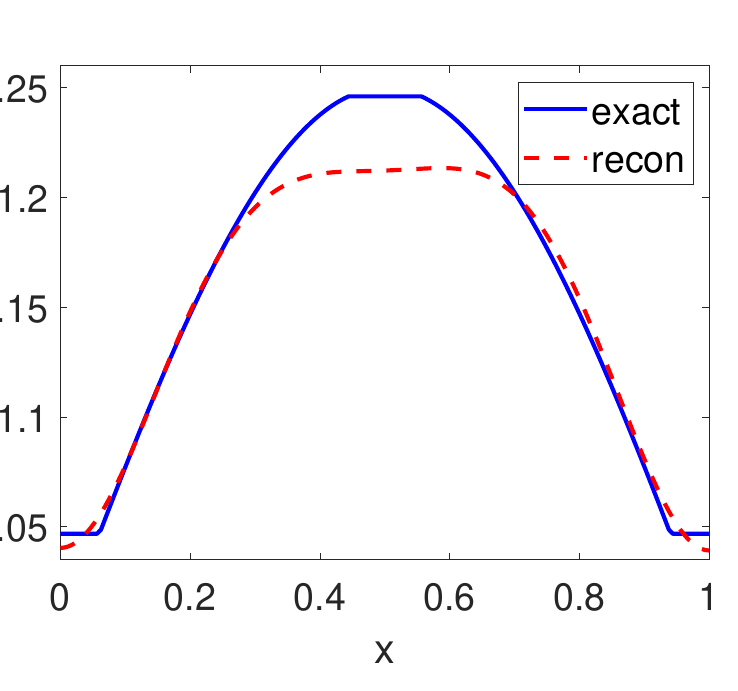} & \includegraphics[width=0.33\textwidth,trim={4.5cm 9cm 4cm 9cm}, clip]{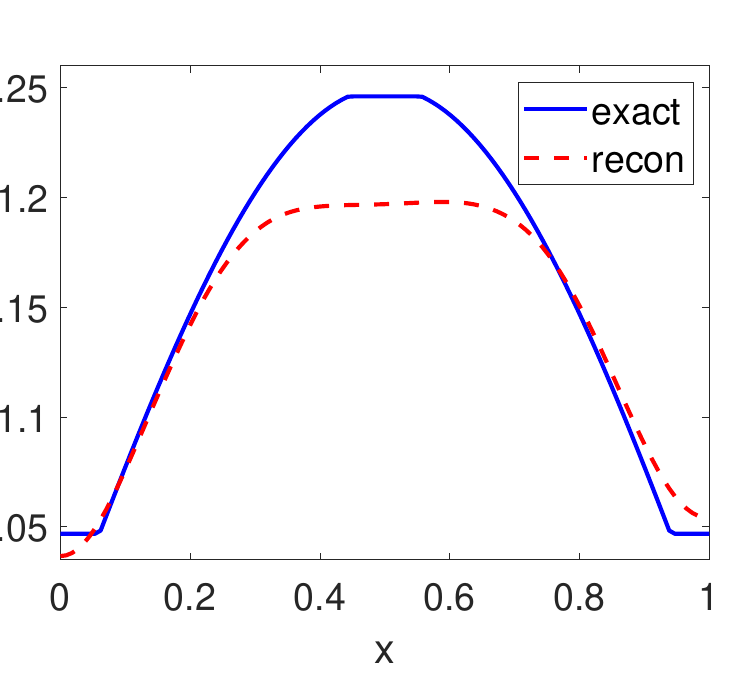}\\
(a) $\epsilon=1.00${\rm e}-3 & (b) $\epsilon=5.00${\rm e}-3 & (c) $\epsilon=1.00${\rm e}-2
\end{tabular}
\caption{The numerical reconstructions for Example \ref{exam:1d}(a) at three noise levels. From top to bottom, the results are for $\alpha=0.25$, $0.50$ and $0.75$.}\label{fig:recon-frac1d}
\end{figure}

Next we give a two-dimensional example.
\begin{example}\label{exam:2d}
Let $\Omega$ be the unit disk centered at the origin $(0,0)$,  $q^\dag(x,y)= \max(\min(1+\frac14\cos\big(\frac{\pi}{2}(x^2+y^2)\big), \frac{319}{256}), \frac{71}{64})$, $u_0(x,y)=1-x^2-y^2$ and $f\equiv1$. Fix $T=2.00$, and $\alpha=0.25$, $0.50$ and $0.75$.
\end{example}

Like before, the exact data is generated using a finer mesh and the discrete optimization problem is solved on coarser meshes. We observe a steady convergence for both $e_q$ and $e_u$: the convergence rate of $e_u$ is slightly slower than the first order; but the convergence rate of $e_q$ is again much higher than theoretical one. Fig \ref{fig:recon-2d} shows exemplary reconstructions and the pointwise error $e:=q^\dag-q_h^*$ at two noise levels $\epsilon=5.19$e-3 and $8.75$e-4.

\begin{table}[htpb]
\centering
\begin{threeparttable}
\caption{Numerical results for Example \ref{exam:2d}, initialized with a mesh with $209$ elements
and the total time step $N=10$.\label{tab:err-2d}}
\begin{tabular}{c|c|cccccc}
\toprule
&$\epsilon$ & 1.00e-2 & 5.19e-3 & 2.14e-3 & 8.75e-4 & rate \\
$\alpha$& $\gamma$   & 1.00e-6 & 2.69e-7 & 4.51e-8 & 7.65e-9 & \\
\midrule
	$0.25$ & $e_q$    & 2.40e-2 & 1.84e-2 & 9.37e-3 & 5.95e-3 & 0.56 \\
	              & $e_u$    & 2.78e-3 & 1.24e-3 & 6.11e-4 & 2.14e-4 & 1.07 \\
	\midrule
	$0.50$ & $e_q$    & 1.97e-2 & 1.43e-2 & 7.58e-3 & 4.52e-3 & 0.59 \\
	              & $e_u$    & 8.31e-4 & 5.50e-4 & 2.02e-4 & 9.97e-5 & 0.85 \\
	\midrule
    $0.75$ & $e_q$    & 2.20e-2 & 1.71e-2 & 9.69e-3 & 4.95e-3 & 0.59 \\
  & $e_u$    & 2.04e-3 & 1.03e-3 & 4.58e-4 & 1.97e-4 & 0.97 \\
    \bottomrule
\end{tabular}
\end{threeparttable}
\end{table}

\begin{figure}[htbp!]
\centering
\setlength{\tabcolsep}{0pt}
\renewcommand{\arraystretch}{0}
\begin{tabular}{ccc}
  \includegraphics[width=0.33\textwidth,trim={4.5cm 9cm 3cm 9cm}, clip]{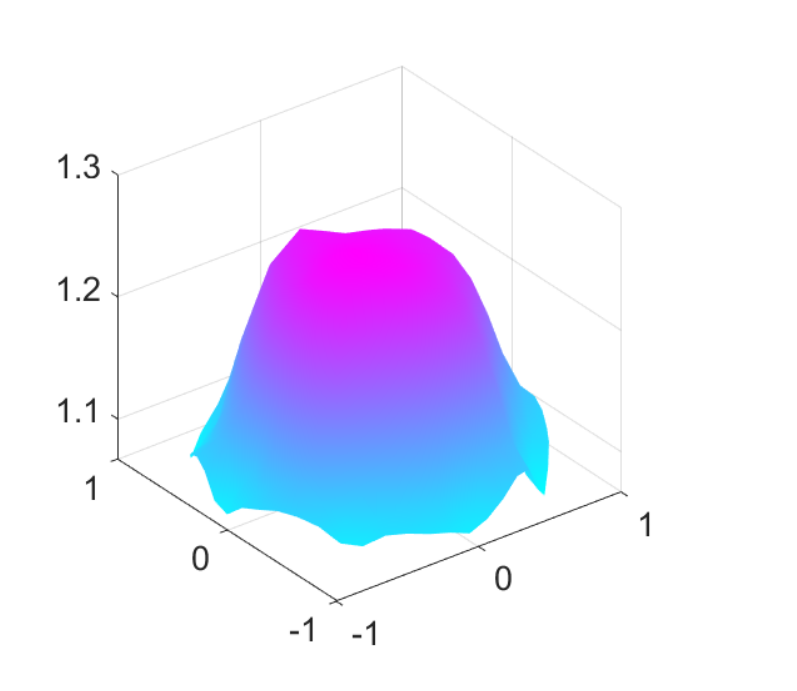} & \includegraphics[width=0.33\textwidth,trim={4.5cm 9cm 3cm 9cm}, clip]{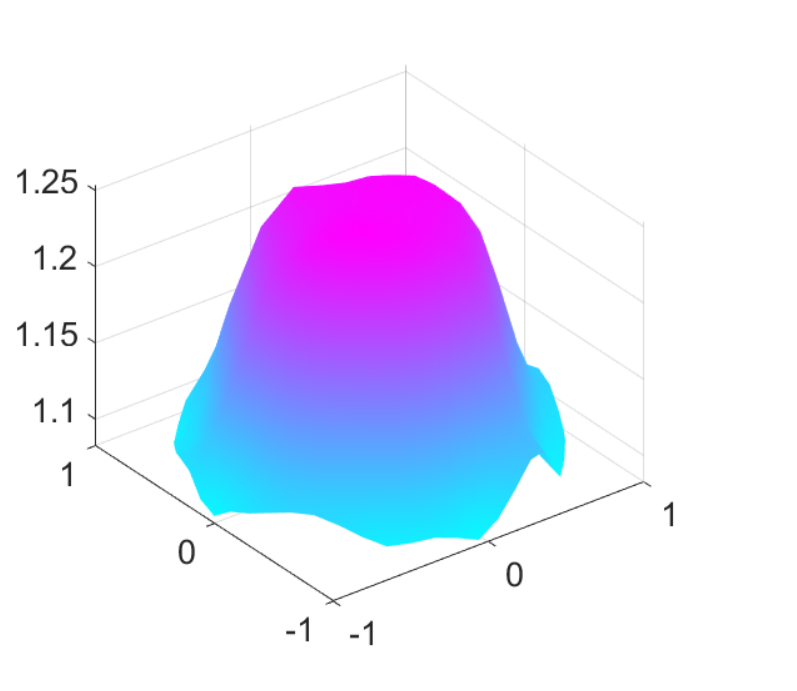} &
  \includegraphics[width=0.33\textwidth,trim={4.5cm 9cm 3cm 9cm}, clip]{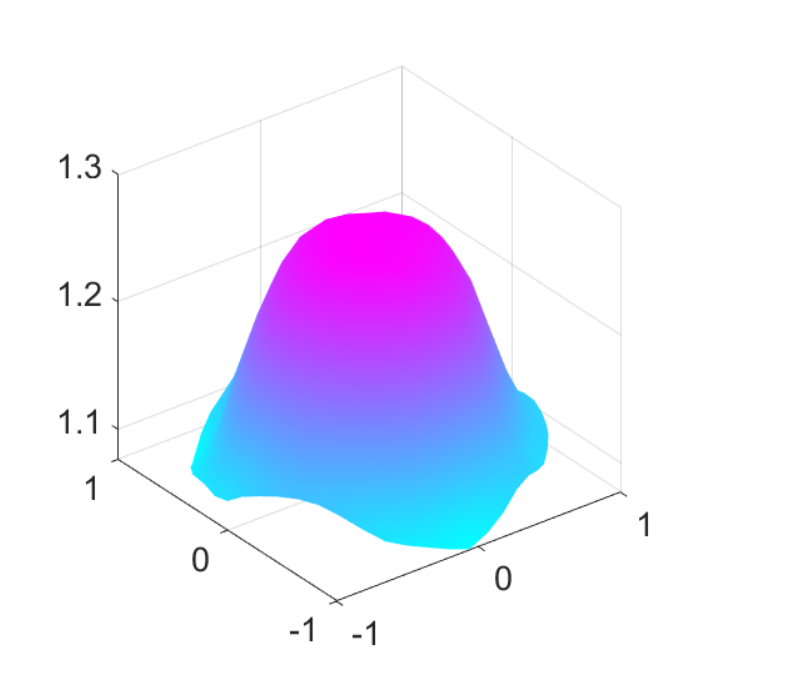} \\
  \includegraphics[width=0.33\textwidth,trim={4.5cm 9cm 3cm 9cm}, clip]{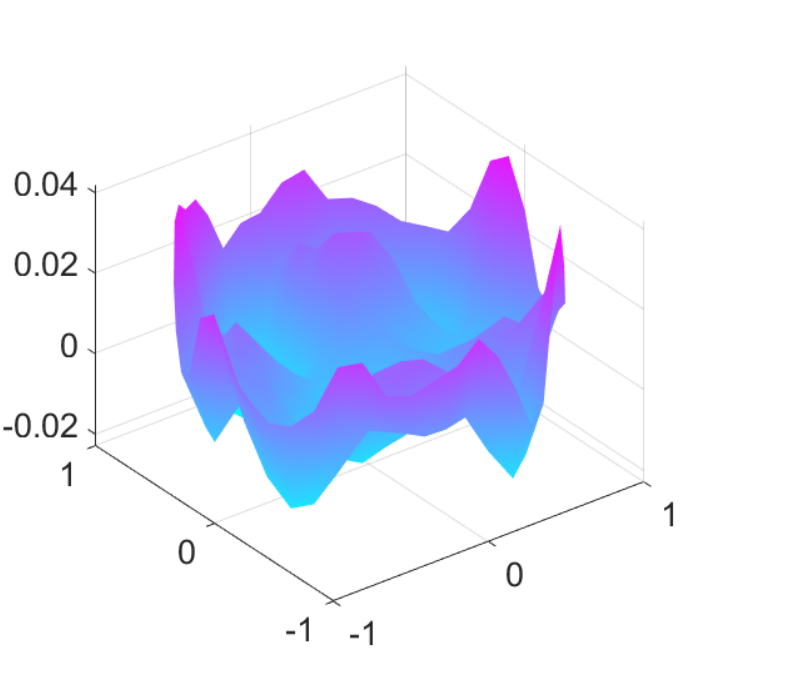} &
  \includegraphics[width=0.33\textwidth,trim={4.5cm 9cm 3cm 9cm}, clip]{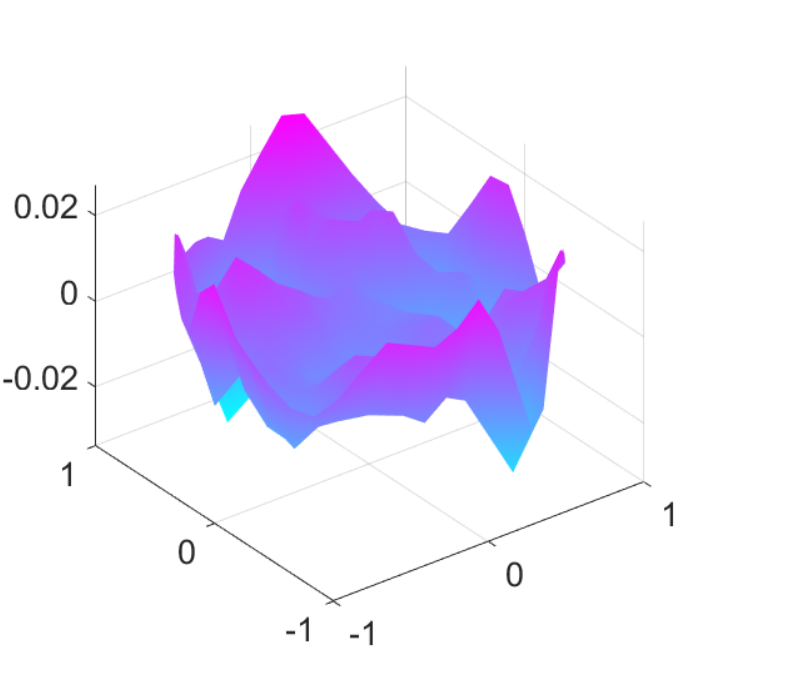} &
  \includegraphics[width=0.33\textwidth,trim={4.5cm 9cm 3cm 9cm}, clip]{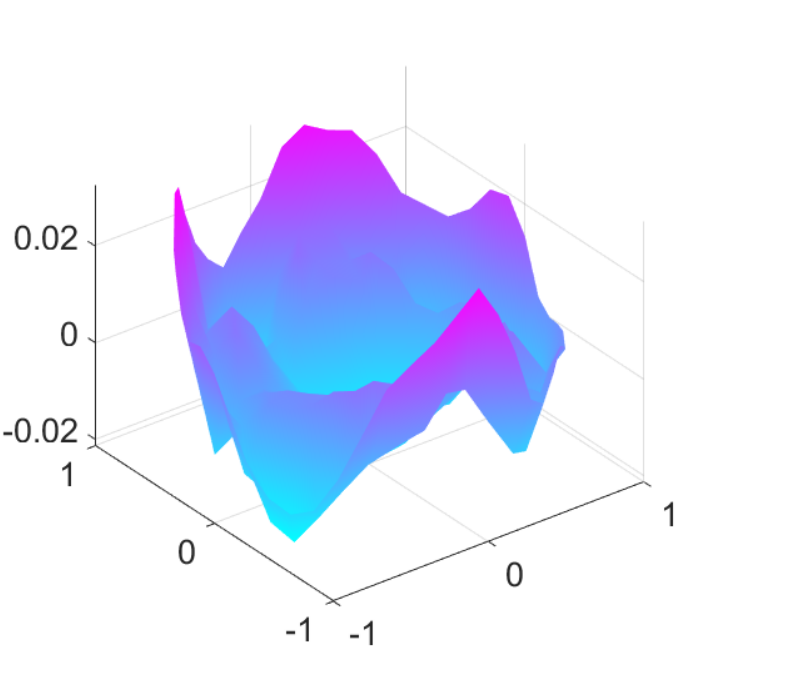} \\
  \includegraphics[width=0.33\textwidth,trim={4.5cm 9cm 3cm 9cm}, clip]{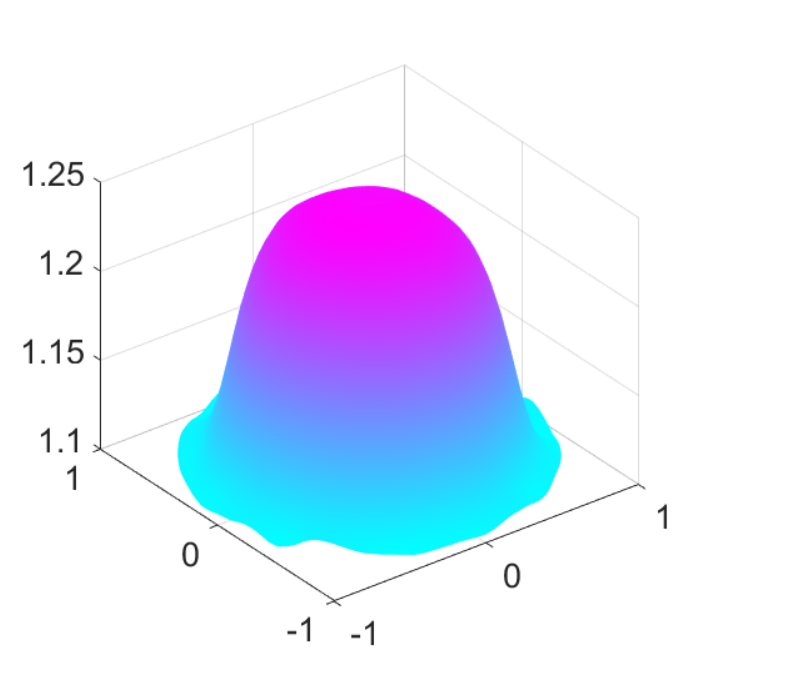} &
  \includegraphics[width=0.33\textwidth,trim={4.5cm 9cm 3cm 9cm}, clip]{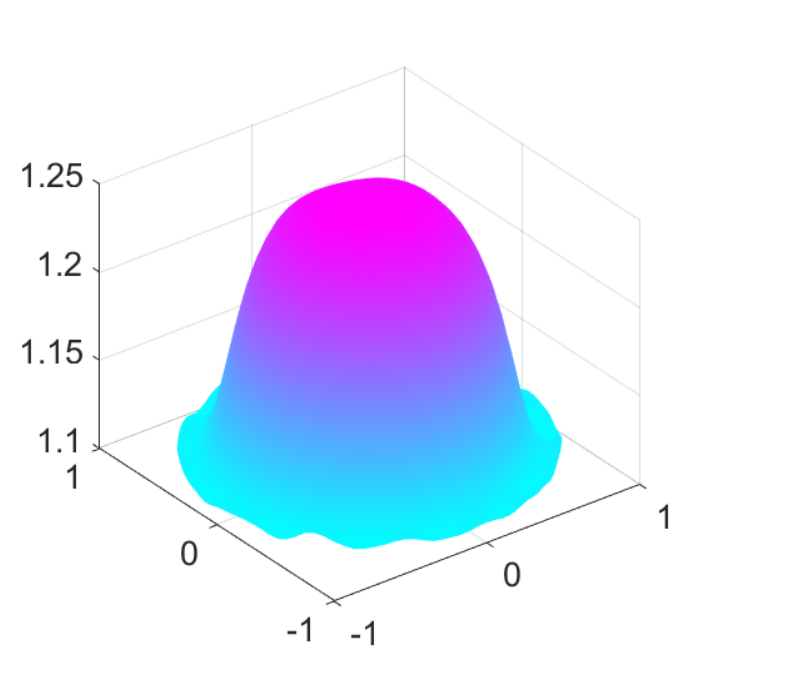} &
  \includegraphics[width=0.33\textwidth,trim={4.5cm 9cm 3cm 9cm}, clip]{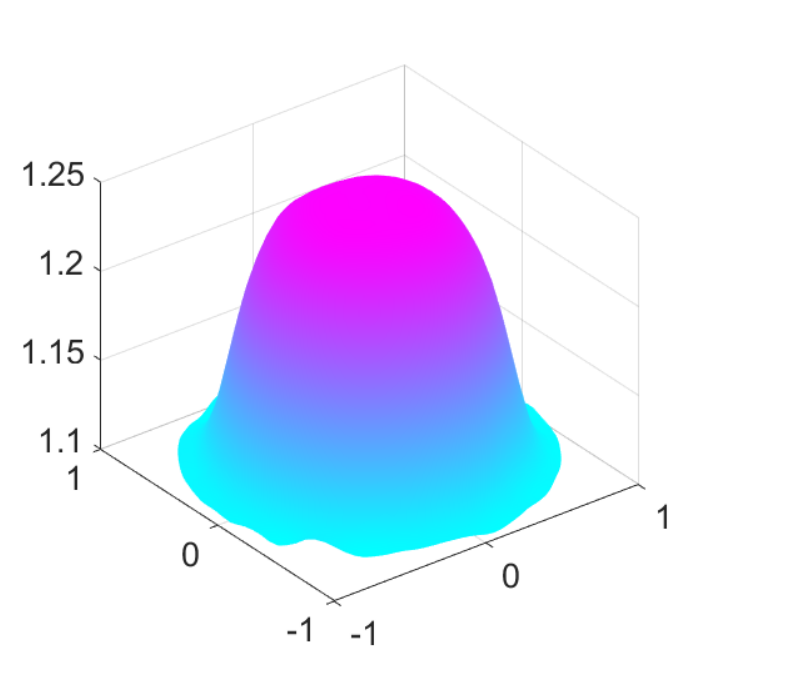}\\
  \includegraphics[width=0.33\textwidth,trim={4.5cm 9cm 3cm 9cm}, clip]{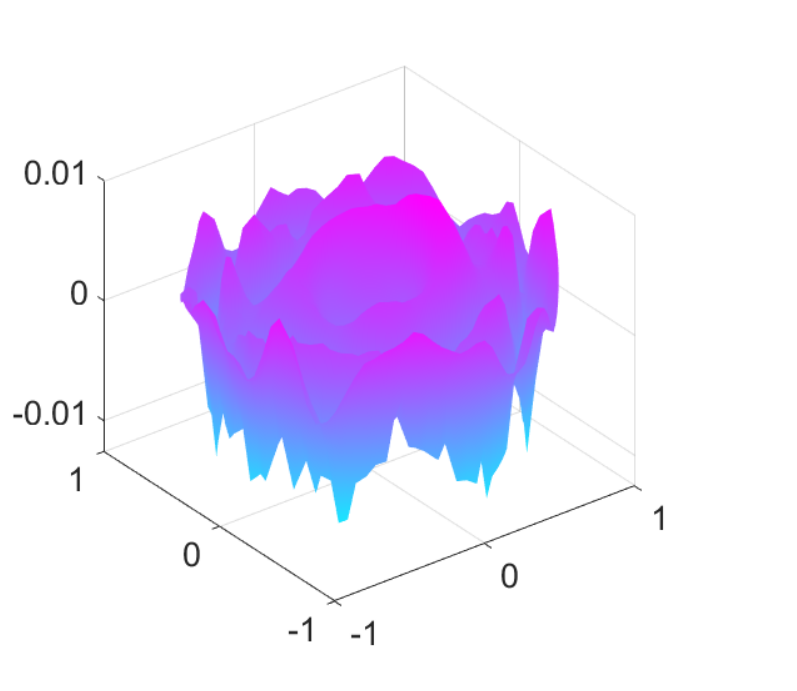} &
  \includegraphics[width=0.33\textwidth,trim={4.5cm 9cm 3cm 9cm}, clip]{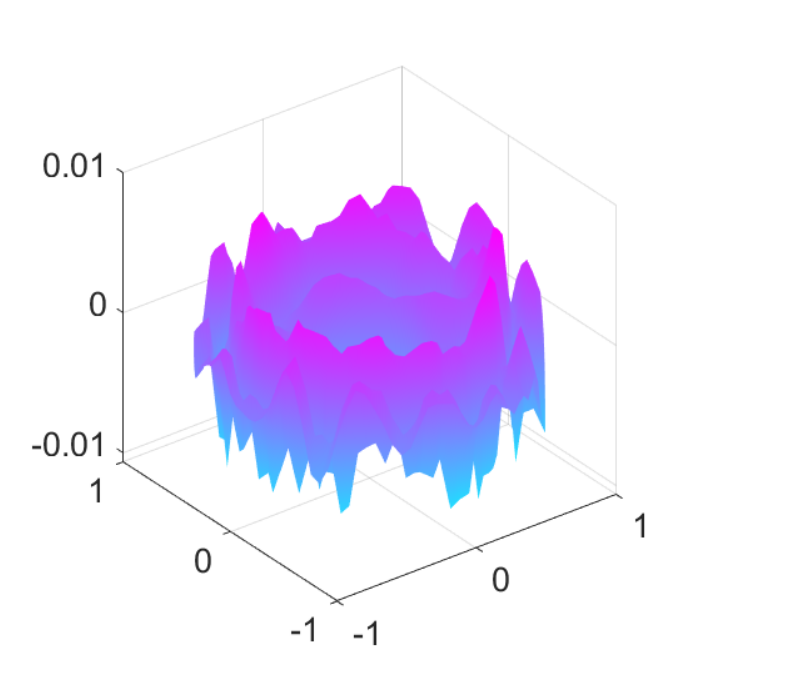} &
  \includegraphics[width=0.33\textwidth,trim={4.5cm 9cm 3cm 9cm}, clip]{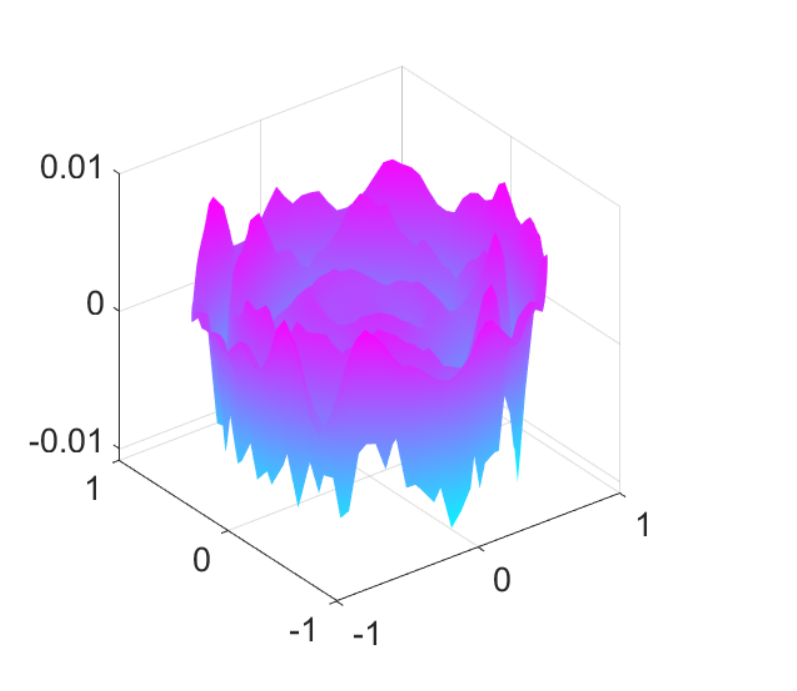}
\end{tabular}
\caption{The numerical reconstructions for Example \ref{exam:2d} with $\epsilon=5.19\text{e-}3$ (top two rows) and $8.75\text{e-}4$ (bottom two rows) and the pointwise error $e:=q^\dag-q_h^*$. From left to right, the results are for $\alpha=0.25$,  $0.50$ and $0.75$.}\label{fig:recon-2d}
\end{figure}

\appendix
\section{Approximation property of Lagrange interpolation}
In this appendix, we prove the approximation property \eqref{inequ:Pi_h approx-2} and \eqref{inequ:Pi_h approx-inf} of the Lagrange nodal interpolation operator $\Pi_h$ in a convex smooth domain.

\begin{lemma}\label{app:lem}
Let $\Omega$ be a convex and smooth domain. Let the polygon $\Omega_h$, the FEM space $V_h$, and the Langrange interpolation operator $\Pi_h: C(\overline\Omega)\rightarrow V_h$ be defined in Section \ref{subsec:Tikhonov regularization problem and its FEM approximation in elliptic system}. Then the estimates
\eqref{inequ:Pi_h approx-2} and \eqref{inequ:Pi_h approx-inf} hold.
\end{lemma}
\begin{proof}
To show the estimate \eqref{inequ:Pi_h approx-inf}, it suffices to show
$$ \|  v - \Pi_h v \|_{L^\infty(\Omega\backslash\Omega_h)}
+
h\|  \nabla(v - \Pi_h v) \|_{L^\infty(\Omega\backslash\Omega_h)}
\le ch \| v \|_{W^{1,\infty}(\Omega)}.$$
By the construction of the space $V_h$, we observe
\begin{align}\label{eqn:nablapih}
\| \nabla \Pi_h v \|_{L^\infty(\Omega\backslash\Omega_h)} \le \| \nabla \Pi_h v \|_{L^\infty(\Omega_h)} \le \| \nabla  v \|_{L^\infty(\Omega_h)} \le \| \nabla  v \|_{L^\infty(\Omega)}.
\end{align}
Moreover, since $\text{dist}(x,\Omega_h) \le ch^2$ for all $x\in\partial\Omega$, from the estimate \eqref{eqn:nablapih}, we derive
\begin{align*}
&\|  v - \Pi_h v \|_{L^\infty(\Omega\backslash\Omega_h)} \le \| v - \Pi_h v  \|_{L^\infty(\partial\Omega_h)} + c h^2 \| \nabla (v - \Pi_h v) \|_{L^\infty(\Omega\backslash\Omega_h)}\\
\le& ch \| v \|_{W^{1,\infty}(\Omega_h)} + c h^2 \| \nabla v\|_{L^\infty(\Omega\backslash\Omega_h)} + c h^2\| \nabla \Pi_h v\| \|_{L^\infty(\Omega\backslash\Omega_h)}\le ch  \| v \|_{W^{1,\infty}(\Omega)}.
\end{align*}
Next, we prove the estimate \eqref{inequ:Pi_h approx-2}. The standard trace lemma \cite[Theorem 1.6.6]{BrennerScott:book2008} implies
\begin{align}\label{eqn:nablapih2}
\|  v - \Pi_h v \|_{L^2(\partial\Omega_h)}
+ h \| \nabla(v - \Pi_h v)\|_{L^2(\partial\Omega_h)}  \le ch^{\frac32} \| v \|_{H^2(\Omega_h)} \le
ch^{\frac32} \| v \|_{H^2(\Omega)}.
\end{align}
Since $\text{dist}(x,\Omega_h) \le ch^2$ for all $x\in\partial\Omega$ and $\Pi_h v$ is piecewise linear, we have
\begin{align*}
\|  v - \Pi_h v \|_{L^2(\Omega\backslash\Omega_h)}^2 &\le ch \| v - \Pi_h v  \|_{L^2(\partial\Omega_h)}^2 + c h^2 \| \nabla (v - \Pi_h v) \|_{L^2(\Omega\backslash\Omega_h)}^2\\
&\le ch \| v - \Pi_h v  \|_{L^2(\partial\Omega_h)}^2
+ c h^3 \| \nabla (v - \Pi_h v) \|_{L^2(\partial\Omega_h)}^2
+ ch^4 \| v \|_{H^2(\Omega\backslash\Omega_h)}^2.
\end{align*}
Then applying \eqref{eqn:nablapih2} gives
$\|  v - \Pi_h v \|_{L^2(\Omega\backslash\Omega_h)} \le ch^2 \| v \|_{H^2(\Omega)}$.
The bound on $\|\nabla(v-\Pi_h v)\|_{L^2(\Omega)}$ follows similarly.
\end{proof}

\section{Conjugate gradient method}
Now we briefly describe the conjugate gradient algorithm \cite{Alifanov:1995}
for minimizing the regularized problem.  The main effort of the algorithm at
each step is to compute the gradient $J'_\gamma(q)$ of the objective $J_\gamma$. This can
be achieved using the adjoint technique. Specifically, let $v(q)$ solve the
modified adjoint equation \cite{cheng2020inverse}
\begin{equation}\label{eqn:adj}
    \left\{\begin{aligned}
	^R\partial_{T+\tau}^\alpha v-\nabla\cdot(q\nabla v) &= (u(q)(T)-z^\delta)\delta_T(t), \ &&\mbox{in}\ \Omega\times(0,T+\tau), \\
     v&=0, \ &&\mbox{on}\ \partial\Omega\times(0,T+\tau), \\
     v(T+\tau)&=0, \ &&\mbox{in}\ \Omega.
   \end{aligned}\right.
\end{equation}
Here $\delta_T(\cdot)$ is Dirac function in $t$ concentrated at $t=T$, and $_t^R\partial_{T}^\alpha v$ is defined by
$_t^R\partial_{T}^\alpha v(t) = -\frac{\rm d}{{\rm d} t}\frac{1}{\Gamma(1-\alpha)}\int_t^T (s-t)^{-\alpha}u(s){\rm d} s$.
Then the $L^2(\Omega)$ gradient $J_\gamma'(q)$ is given by
\begin{equation}\label{eqn:grad}
    J_\gamma'(q) = \nabla u(q)(T)\cdot \nabla v(q)(0) - \gamma\Delta q,
\end{equation}
and the descent direction $g^k=-(-\Delta)^{-1}J_\gamma'(q^k)$. The conjugate gradient direction $d_k$ is given by
\begin{equation}\label{eqn:cg-dir}
    d_k = \beta_k d_{k-1} +g^k, \quad \mbox{with } \beta_k = \|g^{k}\|_{L^2(\Omega)}^2/\|g^{k-1}\|_{L^2(\Omega)}^2,
\end{equation}
with the convention $\beta_0=0$. To select the step size $s$ at Step 6, we employ a simple strategy
by linearizing the direct problem \eqref{eqn:fde} along the direction $d_k$. The operator $P_\mathcal{A}$ at line 7 denotes the pointwise projection into the set $\mathcal{A}$.

\begin{algorithm}[hbt!]
\caption{Conjugate gradient method for problem \eqref{eqn:obj}--\eqref{eqn:weak}.\label{alg:cgm}}
\begin{algorithmic}[1]
\STATE Set the maximum iteration number $K$, and choose $q^0$.
\FOR{$k=1,\ldots,K$}
   \STATE Solve for $u(q^k)$ the solution to problem \eqref{eqn:fde} with $q=q^k$.
   \STATE Solve for $v(q^k)$ the solution to the modified adjoint problem \eqref{eqn:adj} with $q=q^k$;
   \STATE Compute the gradient $J_\gamma'(q^k)$ via \eqref{eqn:grad}, and the descent direction $d_k$ via \eqref{eqn:cg-dir};
   \STATE Compute the step length $s_k$;
   \STATE Update the diffusion coefficient by $q^{k+1}=P_\mathcal{A}(q^k+s_kd_k)$;
   \STATE Check the stopping criterion.
\ENDFOR
\end{algorithmic}
\end{algorithm}

\bibliographystyle{abbrv}
\bibliography{reference}

\begin{thebibliography}{10}

\bibitem{adams1992field}
E.~E. Adams and L.~W. Gelhar.
\newblock Field study of dispersion in a heterogeneous aquifer: 2. spatial
  moments analysis.
\newblock {\em Water Resources Res.}, 28(12):3293--3307, 1992.

\bibitem{Adams2003Sobolev}
R.~A. Adams and J.~J.~F. Fournier.
\newblock {\em {Sobolev Spaces}}.
\newblock Elsevier/Academic Press, Amsterdam, second edition, 2003.

\bibitem{Alessandrini:2020}
G.~Alessandrini.
\newblock A small collection of open problems.
\newblock {\em Rend. Istit. Mat. Univ. Trieste}, 52:591--600, 2020.

\bibitem{AlessandriniVessella:1985}
G.~Alessandrini and S.~Vessella.
\newblock Error estimates in an identification problem for a parabolic
  equation.
\newblock {\em Boll. Un. Mat. Ital. C (6)}, 4(1):183--203, 1985.

\bibitem{Alifanov:1995}
O.~M. Alifanov, E.~A. Artyukhin, and S.~V. Rumyantsev.
\newblock {\em Extreme {M}ethods for {S}olving {I}ll-{P}osed {P}roblems with
  {A}pplications to {I}nverse {H}eat {T}ransfer {P}roblems}.
\newblock Begell House, New York, 1995.

\bibitem{bakaev2001maximum}
N.~Y. Bakaev.
\newblock Maximum norm resolvent estimates for elliptic finite element
  operators.
\newblock {\em BIT}, 41(2):215--239, 2001.

\bibitem{bakaev2003maximum}
N.~Y. Bakaev, V.~Thom\'{e}e, and L.~B. Wahlbin.
\newblock Maximum-norm estimates for resolvents of elliptic finite element
  operators.
\newblock {\em Math. Comp.}, 72(244):1597--1610, 2003.

\bibitem{Bonito:2017}
A.~Bonito, A.~Cohen, R.~DeVore, G.~Petrova, and G.~Welper.
\newblock Diffusion coefficients estimation for elliptic partial differential
  equations.
\newblock {\em SIAM J. Math. Anal.}, 49(2):1570--1592, 2017.

\bibitem{BrennerScott:book2008}
S.~C. Brenner and L.~R. Scott.
\newblock {\em {The Mathematical Theory of Finite Element Methods}}.
\newblock Springer, New York, third edition, 2008.

\bibitem{cheng2009uniqueness}
J.~Cheng, J.~Nakagawa, M.~Yamamoto, and T.~Yamazaki.
\newblock Uniqueness in an inverse problem for a one-dimensional fractional
  diffusion equation.
\newblock {\em Inverse problems}, 25(11):115002, 2009.

\bibitem{cheng2020inverse}
X.~Cheng, L.~Yuan, and K.~Liang.
\newblock Inverse source problem for a distributed-order time fractional
  diffusion equation.
\newblock {\em J. Inverse Ill-Posed Probl.}, 28(1):17--32, 2020.

\bibitem{crouzeix1994resolvent}
M.~Crouzeix, S.~Larsson, and V.~Thom\'{e}e.
\newblock Resolvent estimates for elliptic finite element operators in one
  dimension.
\newblock {\em Math. Comp.}, 63(207):121--140, 1994.

\bibitem{CrouzeixThomee:1987}
M.~Crouzeix and V.~Thom\'{e}e.
\newblock The stability in {$L_p$} and {$W^1_p$} of the {$L_2$}-projection onto
  finite element function spaces.
\newblock {\em Math. Comp.}, 48(178):521--532, 1987.

\bibitem{EggerHofmann:2018}
H.~Egger and B.~Hofmann.
\newblock Tikhonov regularization in {H}ilbert scales under conditional
  stability assumptions.
\newblock {\em Inverse Problems}, 34(11):115015, 17, 2018.

\bibitem{engl1996regularization}
H.~W. Engl, M.~Hanke, and A.~Neubauer.
\newblock {\em {Regularization of Inverse Problems}}.
\newblock Kluwer Academic Publishers Group, Dordrecht, 1996.

\bibitem{fujita1991evolution}
H.~Fujita and T.~Suzuki.
\newblock Evolution problems.
\newblock In {\em {Handbook of Numerical Analysis, {V}ol. {II}}}, Handb. Numer.
  Anal., II, pages 789--928. North-Holland, Amsterdam, 1991.

\bibitem{giona1992fractional}
M.~Giona and H.~E. Roman.
\newblock Fractional diffusion equation for transport phenomena in random
  media.
\newblock {\em Phys. A: Stat. Mech. Appl.}, 185(1-4):87--97, 1992.

\bibitem{hatano1998dispersive}
Y.~Hatano and N.~Hatano.
\newblock Dispersive transport of ions in column experiments: An explanation of
  long-tailed profiles.
\newblock {\em Water Resources Res.}, 34(5):1027--1033, 1998.

\bibitem{Isakov:1991}
V.~Isakov.
\newblock Inverse parabolic problems with the final overdetermination.
\newblock {\em Comm. Pure Appl. Math.}, 44(2):185--209, 1991.

\bibitem{ItoJin:2015}
K.~Ito and B.~Jin.
\newblock {\em Inverse {P}roblems: Tikhonov {T}heory and {A}lgorithms}.
\newblock World Scientific Publishing Co. Pte. Ltd., Hackensack, NJ, 2015.

\bibitem{Jin:2021book}
B.~Jin.
\newblock {\em Fractional {D}ifferential {E}quations---an {A}pproach via
  {F}ractional {D}erivatives}, volume 206 of {\em Applied Mathematical
  Sciences}.
\newblock Springer, Cham, 2021.

\bibitem{JinZhou:2021SINUM}
B.~Jin and Z.~Zhou.
\newblock Error analysis of finite element approximations of diffusion
  coefficient identification for elliptic and parabolic problems.
\newblock {\em SIAM J. Numer. Anal.}, 59(1):119--142, 2021.

\bibitem{JinZhou:2021sicon}
B.~Jin and Z.~Zhou.
\newblock Numerical estimation of a diffusion coefficient in subdiffusion.
\newblock {\em SIAM J. Control Optim.}, 59(2):1466--1496, 2021.

\bibitem{JinZhou:2023book}
B.~Jin and Z.~Zhou.
\newblock {\em {Numerical Treatment and Analysis of Time-Fractional Evolution
  Equations}}.
\newblock Springer, Cham, 2023.

\bibitem{JinZhou:2023IMA}
B.~Jin and Z.~Zhou.
\newblock Recovery of a space-time-dependent diffusion coefficient in
  subdiffusion: stability, approximation and error analysis.
\newblock {\em IMA J. Numer. Anal.}, 43(4):2496--2531, 2023.

\bibitem{keung1998numerical}
Y.~L. Keung and J.~Zou.
\newblock Numerical identifications of parameters in parabolic systems.
\newblock {\em Inverse Problems}, 14(1):83--100, 1998.

\bibitem{kilbas2006theory}
A.~A. Kilbas, H.~M. Srivastava, and J.~J. Trujillo.
\newblock {\em {Theory and Applications of Fractional Differential Equations}}.
\newblock Elsevier Science B.V., Amsterdam, 2006.

\bibitem{LarssonThomee:2003}
S.~Larsson and V.~Thom\'{e}e.
\newblock {\em Partial {D}ifferential {E}quations with {N}umerical {M}ethods}.
\newblock Springer-Verlag, Berlin, 2003.

\bibitem{li2017maximal}
B.~Li and W.~Sun.
\newblock Maximal {$L^p$} analysis of finite element solutions for parabolic
  equations with nonsmooth coefficients in convex polyhedra.
\newblock {\em Math. Comp.}, 86(305):1071--1102, 2017.

\bibitem{li2012numerical}
G.~Li, W.~Gu, and X.~Jia.
\newblock Numerical inversions for space-dependent diffusion coefficient in the
  time fractional diffusion equation.
\newblock {\em J. Inverse Ill-Posed Probl.}, 20(3):339--366, 2012.

\bibitem{li2013simultaneous}
G.~Li, D.~Zhang, X.~Jia, and M.~Yamamoto.
\newblock Simultaneous inversion for the space-dependent diffusion coefficient
  and the fractional order in the time-fractional diffusion equation.
\newblock {\em Inverse Problems}, 29(6):065014, 36, 2013.

\bibitem{LionsMagenes:1972V1}
J.-L. Lions and E.~Magenes.
\newblock {\em {Non-Homogeneous Boundary Value Problems and Applications.
  {V}ol. {I}}}.
\newblock Die Grundlehren der mathematischen Wissenschaften, Band 181.
  Springer-Verlag, New York-Heidelberg, 1972.
\newblock Translated from the French by P. Kenneth.

\bibitem{lubich1986discretized}
C.~Lubich.
\newblock Discretized fractional calculus.
\newblock {\em SIAM J. Math. Anal.}, 17(3):704--719, 1986.

\bibitem{metzler2000random}
R.~Metzler and J.~Klafter.
\newblock The random walk's guide to anomalous diffusion: a fractional dynamics
  approach.
\newblock {\em Phys. Rep.}, 339(1):1--77, 2000.

\bibitem{nigmatullin1986realization}
R.~R. Nigmatullin.
\newblock The realization of the generalized transfer equation in a medium with
  fractal geometry.
\newblock {\em Physica Status Solidi (b)}, 133(1):425--430, 1986.

\bibitem{palencia1996maximum}
C.~Palencia.
\newblock Maximum norm analysis of completely discrete finite element methods
  for parabolic problems.
\newblock {\em SIAM J. Numer. Anal.}, 33(4):1654--1668, 1996.

\bibitem{stewart1974generation}
H.~B. Stewart.
\newblock Generation of analytic semigroups by strongly elliptic operators.
\newblock {\em Trans. Amer. Math. Soc.}, 199:141--162, 1974.

\bibitem{Thome2006GalerkinFE}
V.~Thom\'{e}e.
\newblock {\em Galerkin {F}inite {E}lement {M}ethods for {P}arabolic
  {P}roblems}.
\newblock Springer-Verlag, Berlin, second edition, 2006.

\bibitem{Triki:2021}
F.~Triki.
\newblock Coefficient identification in parabolic equations with final data.
\newblock {\em J. Math. Pures Appl. (9)}, 148:342--359, 2021.

\bibitem{wang2010error}
L.~Wang and J.~Zou.
\newblock Error estimates of finite element methods for parameter
  identifications in elliptic and parabolic systems.
\newblock {\em Discrete Contin. Dyn. Syst. Ser. B}, 14(4):1641--1670, 2010.

\bibitem{yamamoto2001simultaneous}
M.~Yamamoto and J.~Zou.
\newblock Simultaneous reconstruction of the initial temperature and heat
  radiative coefficient.
\newblock {\em Inverse problems}, 17(4):1181, 2001.

\bibitem{zhang2015undetermined}
Z.~Zhang.
\newblock An undetermined coefficient problem for a fractional diffusion
  equation.
\newblock {\em Inverse Problems}, 32(1):015011, 21, 2016.

\bibitem{zhang2022identification}
Z.~Zhang, Z.~Zhang, and Z.~Zhou.
\newblock Identification of potential in diffusion equations from terminal
  observation: analysis and discrete approximation.
\newblock {\em SIAM J. Numer. Anal.}, 60(5):2834--2865, 2022.

\end{thebibliography}

\end{document}